\documentclass[11pt]{amsart}
\usepackage{amsfonts}
\usepackage{amsfonts}
\usepackage{amsfonts}
\usepackage{amsfonts}
\usepackage{amssymb,latexsym}
\usepackage{enumerate}
\newtheorem{theorem}{Theorem}
\newcommand{\R}{\mathbb R}

\newtheorem{remark}[theorem]{Remark}
\newtheorem{proposition}[theorem]{Proposition}
\newtheorem{lemma}[theorem]{Lemma}
\newtheorem{definition}[theorem]{Definition}

\newtheorem{Thm}{Theorem}[section]

\newtheorem{Lem}[Thm]{Lemma}
\newtheorem{Rem}[Thm]{Remark}
\newtheorem{Cor}[Thm]{Corollary}

\numberwithin{equation}{section}
\def\i {{1 \leq i \leq n-1}}
\def\j {{1 \leq j \leq n-1}}
\def\ii {{2 \leq i \leq n-1}}
\def\a {\alpha}
\def\b{\beta}
\def\t{\theta}
\def\n{\nabla}
\def\dfrac{\displaystyle\frac}
\def\dint{\displaystyle\int}

\begin{document}
\setlength{\baselineskip}{1.2\baselineskip}
\title  [Gaussian Curvature estimates for the convex level sets]
{Gaussian Curvature estimates for the convex level sets of solutions for some
nonlinear elliptic partial differential equations}
\author{Pei-He Wang}
\address{School of Mathematical Sciences\\
         Qufu Normal University\\
         Qufu,273165, Shandong Province, China}
\email{peihewang@hotmail.com}
\thanks{Research of the first author was supported by NSF of Shandong No.Q2008A08, research of the second author was supported by NSFC No.10671186 and No.10871187.}
\author{Wei Zhang}
\address{Department of mathematics\\
         University of Science and Technology of China\\
         Hefei,230026, Anhui Province, China.}
         \email{zhwmath@mail.ustc.edu.cn}
\maketitle

\begin{abstract}
We give a lower bound for the Gaussian curvature of convex level sets of minimal  graphs and the solutions to semilinear elliptic equations with the
norm of boundary gradient and the Gaussian curvature of the boundary.
\end{abstract}

\section{Introduction}
\setcounter{equation}{0}

This paper is the continuation of Ma-Ou-Zhang \cite{MOZ09}. In \cite{MOZ09},
they  studied the Gaussian curvature estimates of the convex level sets of
$p$-harmonic function in convex ring in $\mathbb{R}^n$ with homogeneous
Dirichlet boundary conditions. Utilizing the similar technique as in \cite{MOZ09}, in this paper we study the minimal surface equation and some
semilinear elliptic equations.

For minimal surface equation we find a sharp auxiliary function involving the
Gaussian curvature of the convex level sets. It is a harmonic function in
2-dimensional case. In higher dimensions, it is a superharmonic function after
modifying the gradient terms with locally bounded coefficients. Here we use the Laplace-Beltrami operator on the minimal graph (see \cite{GT98}).

For the semilinear elliptic equation with suitable structure conditions, we
find a similar auxiliary function such that it is a superharmonic function in
domain after modifying the gradient terms with locally bounded coefficients.
From these results, we can get the Gaussian curvature estimates of the convex
level sets with the norm of boundary gradient and the Gaussian curvature of the
boundary. For example we obtain the lower bound estimates for the Gaussian
curvature of the level sets of the solutions for a class of similinear elliptic
equations, the strict convexity of its level sets was obtained by
Caffarelli-Spruck \cite{CS82} and Korevaar \cite{Kor}.

The geometry of the level sets of the solutions of elliptic partial
differential equations has been studied for a long time. For instance, Ahlfors
\cite{AH} contains the well-known result that level curves of Green function on
simply connected convex domain in the plane are the convex Jordan curves. In
1931, Gergen \cite{Ge} proved the star-shapeness of the level sets of Green
function on 3-dimensional star-shaped domain. In 1956, Shiffman \cite{Sh56}
studied the minimal surface in $\mathbb{R}^3$. In 1957, Gabriel \cite{Ga57} proved that the level sets of the Green function on a 3-dimensional bounded convex domain are strictly convex. Lewis \cite{Lew77} extended Gabriel's result to $p$-harmonic functions in higher dimensions. Caffarelli-Spruck \cite{CS82}
generalized the Lewis \cite{Lew77} results to a class of semilinear elliptic
partial differential equations. Motivated by the result of Caffarelli-Friedman
\cite{CF85}, Korevaar \cite{Kor} gave a new proof on the results of Gabriel
\cite{Ga57} and Lewis \cite{Lew77} using the following observation: if the
level sets of the $p$-harmonic function is convex with respect to the gradient
direction $\nabla u$, then the rank of the second fundamental form of the level sets is a constant in all domain. A survey of this subject is given by Kawohl
\cite{Ka85}. For more recent related extensions, please see the papers by
Bianchini-Longinetti-Salani \cite{BLS} and Bian-Guan-Ma-Xu \cite{BGMX}.

Now we turn to the question of quantitative results, that is, curvature
estimates of the level sets of the solutions to such elliptic problems. For
2-dimensional harmonic functions and minimal surfaces with convex level curves,
Longinetti \cite{Lo83} and \cite{Lo87} proved that the curvature of the level
sets attains its minimum on the boundary (see also Talenti \cite{T83} for
related results). Recently, Ma-Ou-Zhang \cite{MOZ09} and Chang-Ma-Yang
\cite{CMY} got the Gaussian curvature and principal curvature estimates of the
convex level sets on higher dimensional harmonic functions, and the estimates
give a new approach to get the convexity of the level sets of harmonic
functions. For the other related results see the papers by Rosay-Rudin
\cite{RRu89} and Dolbeault-Monneau \cite{DM}.

Now we state our main theorems.
\begin{Thm}\label{MinimalThm}
Let $\Omega$ be a smooth bounded domain in ${\mathbb{R}}^n (n \geq2)$ and
$u\in C^4(\Omega)\cap C^{2}(\bar{\Omega})$ be the solution of the following
minimal surface equation,
\begin{equation}\label{MCEquation}
{\rm div}\bigg(\frac{\n u}{\sqrt{1+|\n u|^2}}\bigg)=0\qquad {\text{in}}\ \Omega\subset\mathbb{R}^n.
\end{equation}
Assume $|\nabla u|\neq 0$ in $\Omega$. If the level sets of $u$ are strictly
convex with respect to normal $\nabla u$, and let $K$ be the Gaussian
curvature of the level sets. Then we have the following results.
\begin{enumerate}[(i)]
\item For $n=2$, the function $\bigg(\dfrac{|\n u|^2}{1+|\n u|^2}\bigg)^ {-\frac12}K$ attains its minimum and maximum on the boundary $\partial \Omega$, unless it is a constant.

\item For $n \ge 3$, the function $\bigg(\dfrac{|\n u|^2}{1+|\n u|^2}\bigg) ^{\theta}K$ attains its minimum on the boundary $\partial \Omega$ for $\t=-\dfrac 1 2$ or $\t\ge\dfrac {n-3} 2$, unless it is a constant.
\end{enumerate}
\end{Thm}

\begin{Rem}\label{xuanlianmian}
In Theorem~\ref{MinimalThm}, we can choose $\psi(x)=\bigg(\dfrac{|\n u|^2} {1+|\n u|^2}\bigg)^{-\frac12}K(x)$ as our test function. Now we give an
example to explain our choice on $\psi$.

For $r=|x|>2$, let $u(r,\ \t)=\dint_2^r \dfrac{\mathrm{d}s} {\sqrt{s^{2(n-1)}-1}}$ be the $n-$dimensional catenoid. A simple calculation
shows
\begin{align*}
 |\nabla u|^2(x) = \frac{1}{|x|^{2(n-1)}-1},
\end{align*}
and the Gaussian curvature of the level set at $x$ is $$K(x)=|x|^{1-n}.$$
Hence,
\begin{align*}
 \psi(x)=\bigg(\dfrac{|\n u|^2}{1+|\n u|^2}\bigg)^{-\frac 1 2}K(x)\equiv1.
\end{align*}
From the above calculation, one know the choice $\psi(x)= \bigg(\dfrac{|\n u|^2} {1+|\n u|^2}\bigg)^{-\frac12}K(x)$ is sharp. Moreover in 2-dimensional case, under the assumption that $|\nabla u|\neq 0$ in $\Omega$ we shall prove in Corollary~\ref{2-dim} that the function $\bigg(\dfrac{|\n u|^2}{1+|\n u|^2} \bigg)^{-\frac12}K$ is a harmonic function with respect to the Laplace-Beltrami operator on the minimal graph.
\end{Rem}

For the semilinear elliptic equations, under suitable structure conditions on the equation, we have
\begin{Thm}\label{PoissonThm}
Let $\Omega$ be a smooth bounded domain in ${\mathbb{R}}^n(n \geq 2)$ and $u\in C^4(\Omega)\cap C^{2}(\bar{\Omega})$ be a solution of the following equation in $\Omega$, i.e.
\begin{equation}\label{Poisson}
\Delta u =f(x, u) \quad  in \  \Omega,
\end{equation}
where $f\in C^2(\Omega\times\mathbb{R})$ and $f$ is nonnegative. Assume $|\nabla u|\neq 0$ in $\Omega$, and the level sets of $u$ are strictly convex with respect to normal $\nabla u$. Let $K$ be the Gaussian curvature of
the level sets. Then we have the following facts.

(ia) Suppose $f=f(u)$ and $f_u\geq 0$, then the function $|\nabla u|^{-2}K$ attains its minimum on the boundary.

(ib) Suppose $f=f(u)$ and $f_u\le 0$, then the function $|\nabla u|^{n-1}K$ attains its minimum on the boundary.

(ii) Suppose $f=f(x)$ and $t^3f(x)$ is convex with respect to $(x,\,t)\in \Omega \times (0, +\infty)$ (or equivalently $f^{-\frac12}$ is concave for $f$ positive), then the function $|\nabla u|^{n-1}K$ attains its minimum on the boundary.
\end{Thm}
If $u$ is a solution for \eqref{Poisson} with convex level sets with respect to normal $\nabla u$, then we shall prove a useful fact that the norm of gradient $|\n u|$ attains its maximum and minimum on the boundary in Lemma \ref{IncreaseGradident}. Combining this fact and Theorem \ref{PoissonThm} we have the following consequence.
\begin{Cor}\label{CorCaffarelli}
Let $\Omega_0$, $\Omega_1$ be two bounded smooth convex domains in $\R^n(n \ge 2)$ and $\bar\Omega_0\subset\Omega_1$. Let $u$ satisfy
\begin{equation}\label{Quasilinear}
\left\{
\begin{array}{lcl}
              \Delta u = f(u)   &\text{in}&  \Omega=\Omega_0\backslash\bar\Omega_1, \\
                     u = 0   &\text{on}&  \partial \Omega_0,\\
                     u = 1   &\text{on}&  \partial \Omega_1,
\end{array} \right.
\end{equation}
where $f\in C^2([0, 1])$ is nonnegative, non-decreasing and $f(0)=0$. Let $K$ be the Gaussian curvature of the level sets, then we have the following estimate
\begin{align*}
\min_{\Omega} K \geq \bigg(\frac{\min_{\partial \Omega_0}|\n u|}{\max_{\partial \Omega_1}|\n u|}\bigg)^2 \min_{\partial \Omega}K.
\end{align*}
\end{Cor}
In \cite{CS82, Kor}, Caffarelli-Spruck and Korevaar proved the level sets of
solution to \eqref{Quasilinear} are strictly convex with respect to normal
$\nabla u$. In above corollary, we give the quantitative results using the
boundary data.

Assuming $|\nabla u|\neq 0$, Bianchini-Longinetti-Salani \cite{BLS} proved  the convexity of the level sets of solution $u$ for some semilinear elliptic
equation in convex ring with homogeneous Dirichlet boundary conditions. Then
from the constant rank theorem of the second fundamental form of the level sets in \cite{Kor} (or \cite{BGMX}), we know the level sets are strictly convex. For the Poisson equation, our structure condition is the same as theirs. {\it In
summary, for the level sets of the solutions of some class semilinear elliptic
equations, \cite{BLS} gives the convexity, \cite{Kor} (or \cite{BGMX})
guarantees the strict convexity, at last from Theorem \ref{PoissonThm} we can
obtain its lower bound estimates for the Gaussian curvature of the level sets
via the boundary data.}

Now we turn to the minimal surface equation in convex ring with  homogeneous
Dirichlet boundary conditions. Korevaar (see Remark 13 in \cite{Kor}) proved
the strict convexity of the level sets. Using the Theorem \ref{MinimalThm},  we can utilize the same method in proving Corollary \ref{CorCaffarelli} to obtain
the similar lower bound estimates for the Gaussian curvature of the level set
via the boundary data.

\begin{Cor}\label{MSEquationb}
Let $u$ satisfy
\begin{eqnarray}
\left\{
\begin{array}{lcl}
div(\frac{\nabla u}{\sqrt{1+|\nabla u|^2}})=0 \quad & \text{in}&  \Omega=\Omega_0\backslash\bar\Omega_1, \\
              u =  0  & \text{on}&  \partial \Omega_0,\\
              u =  1  & \text{on}&  \partial \Omega_1,
\end{array} \right.
\end{eqnarray}
where $\Omega_0$ and $\Omega_1$ are bounded convex domains in $\R^n, n \ge 3$, $\bar\Omega_1\subset\Omega_0$. Let $K$ be the Gaussian curvature of the level sets, then we have the following estimate
\begin{align*}
\min_{\Omega} K \geq \bigg( \frac{\min_{\partial \Omega_0}|\n u|} {\max_{\partial \Omega_1}|\n u|} \bigg)\frac{\sqrt{1+\min_{\partial \Omega_0} |\nabla u|^2}}{\sqrt{1+\max_{\partial \Omega_1}|\nabla u|^2}} \min_{\partial \Omega}K.
\end{align*}
\end{Cor}

Let $K$ be the Gaussian curvature of the convex level sets. Set
$$\varphi = \log K(x) + \rho(|\nabla u|^2),$$
where the function $\rho$ will be specified later. We shall show the following elliptic differential inequality
\begin{align*}
L(\varphi) \le 0 \,\quad \mod \nabla \varphi \quad \text{in}\quad \Omega,
\end{align*}
where $L$ is the Linearized operator associated with the equation we discussed and here we have suppressed the terms containing the gradient of $\varphi$ with locally bounded coefficients. By applying the strong minimum principle, we then obtain the main results.

In Section 2, we first give brief definitions on the convexity of the level
sets, then obtain the  curvature matrix $(a_{ij})$ of the level sets of a
function, which appeared in \cite{BGMX}. In Section 3, we give some
preliminaries and formal computations for the proof of our theorems. We prove
the Theorem \ref{MinimalThm} in Section 4. In sections 5 we prove the Theorem
\ref{PoissonThm}. Then in the last section we give the proof of
Corollary~\ref{CorCaffarelli}, and we omit the proof of
Corollary~\ref{MSEquationb}. The main technique in the proof of these theorems
consists of rearranging the second and third derivative terms using the
equation and the first derivative condition for $\varphi$. The key idea is the
Pogorelov's method in a priori estimates for fully nonlinear elliptic
equations.

{\bf Acknowledgment:}  The authors would like to thank Prof. X.N. Ma for useful discussions on this subject. The first named author would like to thank the hospitality of the Department of Mathematics of University of Science and Technology of China.

\section{The curvature matrix of level sets}
\setcounter{equation}{0} \setcounter{theorem}{0}

In this section, we shall give the brief definition on the convexity of the
level sets, then introduce the  curvature matrix $(a_{ij})$ of the level sets
of a function, which appeared in \cite{BGMX}. Firstly, we recall some
fundamental notations in classical surface theory. Assume a surface
$\Sigma\subset\R^n$ is given by the graph of a function $v$ in a domain in
$\R^{n-1}$:
$$x_n = v(x'), x'=( x_1,x_2, \cdots, x_{n-1})\in \R^{n-1}.$$
\begin{definition}\label{defin1.1}
We define the graph of function $x_n = v(x')$ is convex with respect to the upward normal $\vec {\nu} = \frac{1}{W}(-v_1,-v_2, \cdots, -v_{n-1},1) $ if the second fundamental form $b_{ij}= \dfrac{v_{ij}}{W}$ of the graph $x_n = v(x')$ is nonnegative definite, where $W=\sqrt{1+|\nabla v|^2}$.
\end{definition}

The principal curvature $\kappa=(\kappa_1, \cdots, \kappa_{n-1})$ of the graph of $v$, being the eigenvalues of the second fundamental form relative to the first fundamental form. We have the following well-known formula.
\begin{lemma}(\cite{CNS85})
The principal curvature of the graph $x_n = v(x')$ with respect to the upward normal $\vec {\nu}$ are the eigenvalues of the symmetric curvature matrix
\begin{equation}\label{1.6}
a_{il} =\frac{1}{W}\bigg\{v_{il} -\frac{v_iv_jv_{jl}}{W(1+W)}
-\frac{v_lv_kv_{ki}}{W(1+W)} + \frac{v_iv_lv_jv_k
v_{jk}}{W^2(1+W)^2}\bigg\},
\end{equation}
where the summation convention over repeated indices is employed.
\end{lemma}
Now we give the definition of the convex level sets of the function $u$. Let $\Omega$ be a domain in $\R^n$  and $u\in C^2(\Omega)$, its level sets can be usually defined in the following sense.
\begin{definition}\label{defin1.2}
Assume $|\nabla u| \neq 0$ in $\Omega$, we define the level set of $u$ passing through the point $x_o \in \Omega$ as $\Sigma^{u(x_o)} =\{x\in \Omega|u(x) =u(x_o)\}$.
\end{definition}

Now we shall locally work near the point $x_o$ where $|\n u(x_o)|\neq 0$. By implicit function theorem, locally the level set $\Sigma^{u(x_o)}$ could be represented as a graph
$$x_n = v(x'), x'=( x_1,x_2, \cdots, x_{n-1})\in \R^{n-1},$$
and $v(x')$ satisfies the following equation
$$u(x_1,x_2, \cdots, x_{n-1}, v(x_1,x_2, \cdots, x_{n-1})) = u(x_o).$$
Then the first fundamental form of the level set is $g_{ij}=\delta_{ij} + \frac{u_iu_j}{u_n^2}, $ and  $ W = (1+|\n v|^2)^{\frac12} = \frac{|\n u|}{|u_n|}$. The upward normal direction of the level set is
\begin{equation}\label{1.11}
\vec{\nu}= \frac{|u_n|}{|\n u|u_n}(u_1,u_2, \cdots, u_{n-1},u_n).
\end{equation}
Let
\begin{equation}\label{1.14}
h_{ij} =u_n^2 u_{ij}+u_{nn}u_{i}u_j-u_nu_ju_{in}-u_nu_iu_{jn},
\end{equation}
then the second fundamental form  of the level set of function $u$ is
$b_{ij}=\frac{v_{ij}}{W}=-\frac{|u_n| h_{ij}}{|\n u|u_n^3}.$
\begin{definition}\label{defin1.3}
For the function $u\in C^2(\Omega)$ we assume $|\nabla u|\neq 0$ in $\Omega$. Without loss of generality we  can let $u_n(x_o)\neq 0$ for $x_o \in \Omega$. We define locally the level set $\Sigma^{u(x_o)} = \{x \in \Omega| u(x)=u(x_o)\}$ is convex respect to the upward normal direction $\vec{\nu}$ if the second fundamental form $b_{ij}$ is nonnegative definite.
\end{definition}
\begin{remark}\label{remark1.5}
If we let $\n u$ be the upward normal of the level set $\Sigma^{u(x_o)}$ at $x_o$, then $u_n(x_o)>0$ by \eqref{1.11}. And from the definition \ref{defin1.3}, if the level set $\Sigma^{u(x_o)}$ is convex with respect to the normal direction $\n u$, then the matrix $(h_{ij}(x_o))$ is nonpositive definite.
\end{remark}

Now we obtain the representation of the curvature matrix $(a_{ij})$ of the level sets of the function $u$ with the derivative of the function $u$,
\begin{equation}\label{1.17}
a_{ij} =\frac{1}{|\n u|u_n^2}\bigg\{-h_{ij} +\frac{u_iu_lh_{jl}}{W(1+W)u_n^2} +\frac{u_ju_lh_{il}}{W(1+W)u_n^2} -\frac{u_iu_ju_ku_l h_{kl}} {W^2(1+W)^2 u_n^4}\bigg\}.
\end{equation}
From now on we denote
\begin{eqnarray}\label{1.18}
B_{ij} = \frac{u_iu_lh_{jl}}{W(1+W)u_n^2} +\frac{u_ju_lh_{il}}{W(1+W)u_n^2}, \quad C_{ij}= \frac{u_iu_ju_ku_l h_{kl}}{W^2(1+W)^2u_n^4},
\end{eqnarray}
and
\begin{eqnarray}\label{1.19}
A_{ij} = -h_{ij} + B_{ij} - C_{ij},
\end{eqnarray}
then the symmetric curvature matrix of the level sets of $u$ could be represented as
\begin{eqnarray}\label{1.20}
a_{ij} =\frac{1}{|\n u|u_n^2}\big[-h_{ij}+ B_{ij} - C_{ij}\big]
=\frac{1}{|\n u|u_n^2}{A_{ij}}.
\end{eqnarray}
With the above notations, we end this section with the following Codazzi's type formula which will be used in the next sections.
\begin{proposition}\label{Prop2.3}(see \cite{BGMX})
Denote $a_{ij,k}=\frac{\partial a_{ij}}{\partial x_k}$ for $1\leq i,j,k\leq n-1$, then at the point where $u_n=|\n u|>0,\, u_i=0$, $a_{ij,k}$ is commutative in ``$i,j,k$", i.e.
\begin{equation*}\label{1.21}
a_{ij,k}=a_{ik,j}.
\end{equation*}
\end{proposition}
\begin{proof}
Direct calculation shows
\begin{equation}\label{1.22}
a_{ij,k}=-u_n^{-1}u_{ijk}+u_n^{-2}(u_{ij}u_{kn}+u_{ik}u_{jn}+u_{jk}u_{in}).
\end{equation}
The right hand side of (\ref{1.22}) is obviously commutative in ``$i,j,k$".
\end{proof}

\section{Preliminaries}
\setcounter{equation}{0} \setcounter{theorem}{0}

In this section, we shall make some preliminary calculation for a general class of elliptic equations. In the following sections we shall work on some special equation, for example minimal graph equation, Poisson equation and semi-linear
elliptic equations.

Let $\Omega$ be bounded smooth domain in $\R^n(n \ge 2)$. Assume $(F^{\alpha
\beta})$ to  be a smooth positive definite function matrix defined in $\Omega$, and $u \in C^4(\Omega)$ be a  solution which satisfies the following equation
$$\sum_{1\leq \a,\b\leq n}F^{\alpha \beta}u_{\alpha \beta}=f(x, u),$$
where $f\in C^2(\Omega\times\mathbb{R})$ is nonnegative, and $F^{\alpha \beta}$ is diagonal at the point $x_o$ where $u_i(x_o)=0 (1 \leq i \leq n-1)$ and
$u_n(x_o)=|\nabla u|>0$.

We assume the level sets of $u$ are strictly convex with respect to the normal
$\n u$, then the curvature matrix $(a_{ij})$ of the level sets is positive
definite in $\Omega$.

Set $$\varphi=\rho(|\nabla u|^2)+ \log K(x),$$ where $K=\det(a_{ij})$ is the
Gaussian curvature of the level sets and $\rho$ is a smooth function defined
on the interval $(0,\ +\infty)$ to be given later. In the following sections,
for suitable choice of $\rho$ we will derive the following elliptic inequality
\begin{equation}\label{3.4}
\sum_{1\leq \alpha, \beta \leq n}F^{\alpha \beta}\varphi_{\alpha\beta} \leq 0 \qquad \mod \nabla \varphi \quad \text{in} \quad \Omega,
\end{equation}
where we modify the terms of $\nabla \varphi$ with locally bounded coefficients.

In order to prove \eqref{3.4} at an arbitrary point $x_o\in\Omega$, as in Caffarelli-Friedman \cite{CF85}, we choose the normal coordinates at $x_o$. By rotating the coordinate system suitably through $T_{x_o}$, we may assume that $u_i(x_o)=0 (1 \leq i \leq n-1)$ and $u_n(x_o)=|\nabla u|>0$. We can further assume that the matrix $(u_{ij}(x_o)) (1\leq i, j\leq n-1)$ is diagonal and
$u_{ii}(x_o) < 0$. We also choose $T_{x_o}$ to vary smoothly with $x_o$. If we can establish $\eqref{3.4}$ at $x_o$ under the above assumptions, then go back to the original coordinate we find that $\eqref{3.4}$ remains valid with new locally bounded coefficients on $\nabla \varphi$ in $\eqref{3.4}$, depending smoothly on the independent variables. Thus it suffices to establish $\eqref{3.4}$ under the above assumptions.

From now on, all the calculation will be done at the fixed point $x_o$.

By taking derivative of $\varphi$, we have
\begin{equation}\label{varphi1}
\varphi_{\alpha}= \sum_{1\leq i, j\leq n-1}a^{ij}a_{ij, \alpha} + \rho '|\nabla u|^2_\alpha.
\end{equation}
It follows that
\begin{equation}\label{MGrad}
\sum_{1\leq i \leq n-1}a^{ii}a_{ii, \alpha} = \varphi_{\alpha}-2\rho' u_n u_{n\alpha}.
\end{equation}
Differentiating equation \eqref{varphi1} once more, we have
\begin{align*}
\varphi_{\alpha \beta}=&\sum_{1\leq i \leq n-1}a^{ii}a_{ii, \alpha \beta} - \sum_{1\leq i,j\leq n-1}a^{ii}a^{jj}a_{ij, \alpha}a_{ij, \beta} +\rho '' |\nabla u|^2_\alpha|\nabla u|^2_\beta +\rho '|\nabla u|^2_{\alpha \beta},
\end{align*}
hence
\begin{align} \label{MFG}
\begin{split}
\sum_{1\leq \alpha, \beta \leq n} F^{\alpha\beta} \varphi_{\alpha \beta} =\ I+II+III+IV,
\end{split}
\end{align}
where
\begin{align*}
I=& \sum_{1\leq i \leq n-1}a^{ii}\sum_{1\leq \alpha,\beta\leq n}F^
{\alpha\beta}a_{ii, \alpha \beta}, \hspace{1.5cm} II= - \sum_{1\leq
i,j\leq n-1} \sum_ {1\leq \alpha,\beta \leq
n}F^{\alpha\beta}a^{ii}a^{jj}a_{ij,
\alpha}a_{ij, \beta}, \\
III=&~\rho ''\sum_{1\leq\alpha, \beta\leq n}F^{\alpha\beta}|\nabla
u|^2_\alpha|\nabla u|^2_\beta,\hspace{1.65cm} IV=\ \rho
'\sum_{1\leq\alpha, \beta\leq n}F^{\alpha\beta}|\nabla u|^2_{\alpha
\beta}.
\end{align*}
In the rest of this section, we will deal with the four terms above respectively.

For the term $III$, we have
\begin{align}\label{III}
\begin{split}
III=&~\rho ''\sum_{1\leq\alpha, \beta\leq n} F^{\alpha\beta}|\nabla
u|^2_\alpha|\nabla u|^2_\beta=4\rho '' u_n^{2}\sum_{1\leq \alpha,\b
\leq n} F^{\alpha\beta}u_{n\alpha}u_{n\beta}.
\end{split}
\end{align}

In a similar way, for the term $IV$ we obtain
\begin{align}\label{IV}
\begin{split}
IV=&\ \rho '\sum_{1\leq\alpha, \beta\leq n} F^{\alpha\beta} |\nabla
u|^2_{\alpha \beta} \quad = \rho '\sum_{1\leq\alpha, \beta\leq n}
F^{\alpha\beta} \bigg(2u_nu_{n\alpha\beta} +2\sum_{1\leq \gamma \leq n} u_{\gamma\alpha}u_{\gamma\beta}\bigg)\\
=&\ 2\rho ' u_n\sum_{1\leq \a, \b \leq n} F^{\a\b}u_{\a\b n} + 2\rho' \sum_{1 \leq \a,\b \leq n}F^{\a\b}u_{n\alpha}u_{n\beta}+2\rho '\sum_\i \sum_{1\leq \a,\b \leq n} F^{\a\b} u_{\alpha i}u_{\beta i}.
\end{split}
\end{align}

Next, we deal with the term $I$. By \eqref{1.20}, one has
\begin{align}\label{MAa}
A_{ii}=a_{ii}E, \qquad{\text{where}} \ \  E=|\nabla u|u_{n}^2.
\end{align}
Taking the second derivative of \eqref{MAa}, we get
\begin{align*}
A_{ii,\a\b}=~a_{ii,\a\b}E+ a_{ii,\a}E_\b +a_{ii,\b}E_\a +a_{ii}E_{\a\b},
\end{align*}
so
\begin{align}\label{MI}
\begin{split}
I= \sum_{1\leq i \leq n-1}a^{ii}\sum_{1\leq \alpha,\beta\leq n}F^{\alpha\beta} a_{ii, \alpha \beta} =~I_1+I_2+I_3,
\end{split}
\end{align}
where
\begin{align*}
I_1=&~u_n^{-3}\sum_\i a^{ii}\bigg(\sum_{1\leq \a,\b \leq n} F^{\a\b} A_{ii,\a\b}\bigg), \hspace{0.55cm} I_2=-(n-1) u_n^{-3} \sum_{1\leq \a,\b \leq n}F^{\a\b} E_{\a\b}, \\
I_3=&-2u_n^{-3} \sum_{1\leq \a,\b \leq n} F^{\a\b} \bigg(\sum_\i a^{ii} a_{ii,\b} \bigg)E_{\a}.
\end{align*}
Since
\begin{align*}
E_\a= 3u_n^2u_{n\a},\ \ \ E_{\a\b}= 5u_nu_{n\a}u_{n\b} +3u_n^2u_{n\a\b} +u_n\sum_{1\leq \gamma\leq n}u_{\gamma \a} u_{\gamma \b},
\end{align*}
by \eqref{MGrad}, we obtain
\begin{align}\label{MI3}
\begin{split}
I_3=~12\rho '\sum_{1 \leq \a,\b \leq n}F^{\a\b}u_{n\a}u_{n\b}- 6u_n^{-1}\sum_{1 \leq \a,\b \leq n}F^{\a\b} u_{n\alpha}\varphi_\b;
\end{split}
\end{align}
also
\begin{align}\label{MI2}
\begin{split}
I_2 =&-(n-1)u_n^{-3} \bigg(5u_n\sum_{1\leq\a,\b\leq n}
F^{\a\b}u_{n\a} u_{n\b} +3u_n^2 \sum_{1\leq\a,\b\leq n}
F^{\a\b}u_{n\a\b} +u_n \sum_{1\leq\a,\b,\gamma\leq n}
F^{\a\b}u_{\gamma \a}u_{\gamma \b}\bigg)\\
=&-6(n-1)u_n^{-2}\sum_{1\leq\a,\b\leq n} F^{\a\b}u_{n\a} u_{n\b}
-3(n-1)u_n^{-1} \sum_{1\leq\a,\b\leq n} F^{\a\b}u_{\a\b n}\\
&-(n-1)u_n^{-2} \sum_\i \sum_{1\leq \a,\b \leq n} F^{\a\b} u_{\alpha
i}u_{\beta i}.
\end{split}
\end{align}

For the term $I_1$, recalling the definition of $A_{ij}$, i.e. \eqref{1.18} and \eqref{1.19}, at $x_o$ we have
\begin{align*}
C_{ii,\alpha\beta}=0,
\end{align*}
therefore
$$A_{ii,\a\b} = -h_{ii,\a\b} + B_{ii,\a\b}.$$
By \eqref{1.18},
\begin{eqnarray*}
\begin{split}
&\ u_n^{-3}\sum_\i a^{ii}\bigg(\sum_{1\leq \a,\b \leq n}F^{\a\b} B_{ii, \a\b} \bigg)\\
=&\ u_n^{-3}\sum_\i a^{ii}\sum_{1\leq \a,\b \leq n}F^{\a\b} \left(\frac{2\sum_{1\leq l\leq n-1}u_iu_lh_{il}}{W(1+W)u_n^2} \right)_{\a\b}\\
=& -2u_n^{-2}\sum_\i \sum_{1\leq \a,\b \leq n} F^{\a\b} u_{\alpha i} u_{\beta i}.
\end{split}
\end{eqnarray*}

It follows that
\begin{align}\label{MI1}
\begin{split}
I_1=&\ u_n^{-3}\sum_\i a^{ii}\bigg(\sum_{1\leq \a,\b \leq n} F^{\a\b} A_{ii,\a\b} \bigg)\\
=&\ u_n^{-3}\sum_\i a^{ii}\sum_{1\leq \a,\b \leq n} F^{\a\b} (-h_{ii,\a\b}) -2u_n^{-2}\sum_\i \sum_{1\leq \a,\b \leq n} F^{\a\b} u_{\alpha i}u_{\beta i}.
\end{split}
\end{align}

Combining \eqref{MFG}-\eqref{IV}, \eqref{MI}--\eqref{MI1}, it yields
\begin{align} \label{MFG2}
\begin{split}
\sum_{1\leq \alpha, \beta \leq n} F^{\alpha\beta} \varphi_{\alpha \beta} = &~u_n^{-3}\sum_{1\leq i \leq n-1}a^{ii}\sum_{1\leq \alpha, \beta \leq n}F^ {\alpha\beta}(-h_{ii, \alpha \beta}) - \sum_{1\leq i,j\leq n-1} \sum_ {1\leq \alpha,\beta \leq n}F^{\alpha\beta}a^{ii}a^{jj}a_{ij, \alpha} a_{ij, \beta} \\
&+[2\rho 'u_n-3(n-1)u_n^{-1}] \sum_{1\leq \a, \b \leq n} F^{\a\b}u_{n\a\b}\\
 &+[4\rho '' u_n^{2}+14\rho '-6(n-1)u_n^{-2}]\sum_{1\leq \alpha,\b \leq n}
F^{\alpha\beta}u_{n\alpha}u_{n\beta}\\
&+[2\rho '-(n+1)u_n^{-2}]\sum_{1 \leq \a,\b \leq n}\sum_\i F^{\a\b} u_{\alpha i} u_{\beta i} -6u_n^{-1}\sum_{1 \leq \a,\b \leq n}F^{\a\b} u_{n\alpha} \varphi_\b.
\end{split}
\end{align}

Next, we will compute the term
$$u_n^{-3}\sum_{1\leq i \leq n-1}a^{ii}\sum_{1\leq \alpha, \beta \leq n}F^ {\alpha\beta} (-h_{ii, \alpha \beta}).$$
By differentiating \eqref{1.14} twice, we have
\begin{align*}
-h_{ii, \a}=& -u_n^2u_{ii\a}
-2u_nu_{n\a}u_{ii}-u_{nn\a}u_i^2-2u_iu_{i\a}u_{nn}\\
&+2u_nu_{i}u_{ni\a} +2u_{n}u_{i\a}u_{ni}+2u_{n\a}u_{i}u_{ni},
\end{align*}
and
\begin{align*}
-h_{ii, \a\b}=& -u_n^2u_{ii\a\b} -2u_nu_{n\a}u_{ii\b} - 2u_nu_{n\b}
u_{ii\a} +2u_nu_{i\a}u_{ni\b}\\&+2u_nu_{i\b}u_{ni\a}
+2u_nu_{i\a\b}u_{ni} -2u_nu_{n\a\b}u_{ii} +2u_{n\a}u_{i\b}u_{ni}\\&
+ 2u_{n\b}u_{i\a}u_{ni}-2u_{nn}u_{i\b}u_{i\a}-2u_{n\b}u_{n\a}u_{ii},
\end{align*}
since at $x_o$, $u_{ii}=-u_na_{ii}$, therefore
\begin{align}\label{MFh}
\begin{split}
&u_n^{-3}\sum_{1\leq i \leq n-1}a^{ii}\sum_{1\leq\a,\b\leq
n}F^{\a\b} (-h_{ii,\a\b}) \\=& -u_n^{-1}\sum_{1\leq i \leq
n-1}a^{ii}\sum_{1\leq\a, \b\leq n} F^{\a\b}u_{ii\a\b}
-4u_n^{-2}\sum_{1\leq i \leq n-1}a^{ii}\sum_{1\leq\a,\b\leq
n}F^{\a\b}u_{n\a}u_{ii\b} \\&~+ 4u_n^{-2}\sum_{1\leq i \leq
n-1}a^{ii}\sum_{1\leq\a,\b\leq n}F^{\a\b}u_{i\a}u_{ni\b}
+2u_n^{-2}\sum_{1\leq i \leq n-1}a^{ii}u_{ni}\sum_{1\leq\a,\b\leq
n}F^{\a\b}u_{\a\b i}\\&~ +2(n-1)u_n^{-1}\sum_{1\leq\a,\b\leq
n}F^{\a\b}u_{\a\b n} +4u_n^{-3}\sum_{1\leq i \leq
n-1}a^{ii}u_{ni}\sum_{1\leq\a,\b\leq n}F^{\a\b}u_{n\a}u_{i\b}\\
&~-2u_n^{-3}u_{nn}\sum_{1\leq i \leq n-1}a^{ii}\sum_{1\leq\a,\b\leq
n}F^{\a\b}u_{i\b}u_{i\a} +2(n-1)u_n^{-2}\sum_{1\leq\a,\b\leq
n}F^{\a\b}u_{n\b}u_{n\a}.
\end{split}
\end{align}

By \eqref{1.14}, \eqref{1.18} and \eqref{1.19}, we have
\begin{align*}
A_{ii,\a}=-h_{ii,\a} =-2u_nu_{n\a}u_{ii}-u_n^2u_{ii\a}+2u_nu_{i\a}u_{in},
\end{align*}
on the other hand, by \eqref{1.20},
\begin{align*}
A_{ii,\a}=(a_{ii}|\nabla u|u_n^2)_\a=u_n^3a_{ii,\a} +3u_n^2u_{n\a}a_{ii},
\end{align*}
 therefore
\begin{align}\label{uiia}
u_{ii\a}=-u_na_{ii,\a}+2u_n^{-1}u_{ni}u_{i\a} -u_{n\a}a_{ii}.
\end{align}
By \eqref{MGrad} and \eqref{uiia}, we have
\begin{align}\label{uiiaa}
\begin{split}
&-4u_n^{-2}\sum_{\i}a^{ii}\sum_{1\leq\a,\b\leq n}F^{\a\b}u_{n\a}u_{ii\b}\\
=&-4u_n^{-2}\sum_{1\leq\a,\b\leq n}F^{\a\b}u_{n\a}\sum_{\i}a^{ii}(-u_na_{ii,\b} +2u_n^{-1}u_{ni}u_{i\b} +u_n^{-1}u_{n\b}u_{ii})\\
=&~[-8\rho'+4(n-1)u_n^{-2}]\sum_{1\leq \alpha,\b \leq n} F^{\alpha\beta} u_{n\alpha} u_{n\beta} -8u_n^{-3}\sum_{\i}a^{ii}u_{ni} \sum_{1\leq\a,\b\leq n}F^{\a\b}u_{n\a}u_{i\b}\\&+4u_n^{-1}\sum_{1\leq\a,\b\leq n}F^{\a\b} u_{n\a} \varphi_{\b}.
\end{split}
\end{align}
Noticing that $(F^{\a\b})$ is diagonal at the considered point $x_o$, in a similar way we can obtain
\begin{align}\label{uiia2}
\begin{split}
&~4u_n^{-2}\sum_{\i}a^{ii}\sum_{1\leq\a,\b\leq n}F^{\a\b} u_{i\a}u_{ni\b}\\
=&~4u_n^{-2}\sum_{\i}a^{ii}\sum_{1\leq \b \leq n}F^{i\b}u_{ii}u_{ni\b} +4u_n^{-2}\sum_{\i}a^{ii}\sum_{1\leq\b\leq n}F^{n\b}u_{ni}u_{ni\b}\\
=&-4u_n^{-1}\sum_{1\le i\le n-1}F^{ii}u_{iin} +4u_n^{-2}\sum_{\i}a^{ii}u_{ni} \sum_{1\leq\a,\b\leq n}F^{\a\b}u_{\a\b i}\\
&-4u_n^{-2}\sum_{1\leq i, j\leq n-1}a^{ii}u_{ni}F^{jj}u_{jji}.
\end{split}
\end{align}

Combining \eqref{MFG2}-\eqref{MFh} and \eqref{uiiaa}-\eqref{uiia2}, we obtain
\begin{align} \label{J2+J3}
\begin{split}
\sum_{1\leq \alpha, \beta \leq n} F^{\alpha\beta} \varphi_{\alpha
\beta}=&-u_n^{-1}\sum_{\i}a^{ii}\sum_{1\leq \alpha,\beta\leq
n}F^ {\alpha\beta}u_{ii, \alpha \beta} - \sum_{1\leq i,j\leq n-1}
\sum_ {1\leq \alpha,\beta \leq n}F^{\alpha\beta}a^{ii}a^{jj}a_{ij,
\alpha}a_{ij, \beta},\\
&+6u_n^{-2}\sum_{\i}a^{ii}u_{ni}\sum_{1\leq\a,\b\leq n}F^{\a\b}u_{\a\b i} +[2\rho 'u_n-(n-1)u_n^{-1}] \sum_{1\leq \a, \b \leq n} F^{\a\b}u_{\a\b n}\\
&-4u_n^{-1}\sum_{\i}F^{ii}u_{iin}-4u_n^{-2}\sum_{1\leq i,j\leq n-1}a^{ii}u_{ni}F^{jj}u_{jji}\\
&-4u_n^{-3}\sum_{\i}a^{ii}u_{ni}\sum_{1\leq\a,\b\leq n}F^{\a\b}u_{n\a}u_{i\b} -2u_{nn}u_n^{-3}\sum_{\i}a^{ii}\sum_{1\leq\a,\b\leq n}F^{\a\b} u_{i\b} u_{i\a}\\
&+ [2\rho '-(n+1)u_n^{-2}]\sum_{1 \leq \a,\b \leq n}\sum_\i F^{\a\b} u_{\alpha i}u_{\beta i} +[4\rho''u_n^2+6\rho']\sum_{1\leq \alpha,\b \leq n} F^{\alpha\beta} u_{n\alpha}u_{n\beta}\\
&-2u_n^{-1}\sum_{1 \leq \a,\b \leq n}F^{\a\b} u_{n\alpha}\varphi_\b.
\end{split}
\end{align}

Noticed that $F^{nn}u_{nn}=\sum\limits_{1 \leq \a,\b \leq n}F^{\a\b} u_{\a\b} -\sum\limits_{1 \leq j \leq n-1}F^{jj}u_{jj}$, the terms in the fourth line of the formula \eqref{J2+J3} can be computed as
\begin{align}\label{twoterm}
\begin{split}
&-4u_n^{-3}\sum_{\i}a^{ii}u_{ni}\sum_{1\leq\a,\b\leq n} F^{\a\b} u_{n\a} u_{i\b}-2u_{nn}u_n^{-3}\sum_{\i}a^{ii}\sum_{1\leq\a,\b\leq n}F^{\a\b} u_{i\b} u_{i\a}\\
=&-6u_n^{-3}F^{nn}u_{nn}\sum_{\i}a^{ii}u_{ni}^2 +4u_n^{-2}\sum_{1 \leq j \leq n-1}F^{jj}u_{nj}^2 +2u_n^{-2}u_{nn}\sum_{1 \leq j \leq n-1}F^{jj} u_{jj}\\
=&~\bigg(-6u_n^{-3}\sum_{1\leq\a,\b\leq n}F^{\a\b} u_{\a\b} +6u_n^{-3} \sum_{1\leq j \leq n-1}F^{jj} u_{jj}\bigg )\sum_{\i}a^{ii}u_{ni}^2\\
&+4u_n^{-2}\sum_{1 \leq j \leq n-1}F^{jj}u_{nj}^2 +2u_n^{-2}u_{nn} \sum_{1 \leq j \leq n-1}F^{jj}u_{jj}
\end{split}
\end{align}

By inserting \eqref{twoterm} into \eqref{J2+J3}, we can deduce the following formula
\begin{align}\label{H1234}
\sum_{1\leq \alpha, \beta \leq n} F^{\alpha\beta} \varphi_{\alpha \beta} =L_1+L_2+L_3+L_4,
\end{align}
where
\begin{align*}
L_1=&-{u_n}^{-1}\sum_{\i}a^{ii}\sum_{1\leq \alpha,\beta\leq n}F^
{\alpha\beta}u_{\alpha \beta ii} +6u_n^{-2}\sum_{\i} a^{ii}u_{ni}
\sum_{1\leq\a,\b\leq n}F^{\a\b}u_{\a\b i}\\
&-6u_n^{-3}\sum_{1\leq \alpha,\beta\leq n}F^ {\alpha\beta}u_{ \alpha
\beta}\sum_{\i}a^{ii}u_{ni}^2+[2\rho 'u_n-(n-1)u_n^{-1}] \sum_{1\leq
\a, \b \leq n} F^{\a\b}u_{n\a\b}\,;
\end{align*}
\begin{align*}
L_2=&- \sum_{1\leq i,j,k\leq n-1} F^{kk}a^{ii}a^{jj}a_{ij, k}^2-F^{nn} \sum_{1\leq i,j\leq n-1} a^{ii}a^{jj}a_{ij, n}^2
-4u_n^{-1}\sum_{\i}F^{ii}u_{iin}\\
&-4u_n^{-2}\sum_{1\leq i,j\leq n-1}a^{ii}u_{ni}F^{jj}u_{jj i};\\
L_3=&~[4\rho''u_n^2+6\rho']F^{nn}u_{nn}^2+ [4\rho''u_n^2+6\rho'+4u_n^{-2}] \sum_{\i}F^{ii}u_{ni}^2\\
&+6u_n^{-3}\sum_{\j}F^{jj}u_{jj}\sum_{\i}a^{ii}u_{ni}^2+2u_n^{-2}u_{nn}\sum_{\j}F^{jj}u_{jj}\\
&+[2\rho '-(n+1)u_n^{-2}]F^{nn}\sum_{\i}u_{ni}^2
+[2\rho'-(n+1)u_n^{-2}]\sum_{\j}F^{jj}u_{jj}^2;\\
L_4=&-2u_n^{-1}\sum_{1\leq\a,\b\leq n}F^{\a\b}u_{n\a}\varphi_{\b}.
\end{align*}

Let us state the following lemma.
\begin{lemma}\label{Form2}
Let $A, B, C, D$ be four constants  and $A>0, C>0$. Denote
\begin{align}
\begin{split}
M_1=& -A\sum_{1\leq i,j,k\leq n-1}a^{ii}a^{jj}a_{ij, k}^2
+4B\sum_{1\leq i,j\leq n-1}a^{jj}u_{nj}a_{ii, j},\\
M_2=&-C\sum_{1\leq i,j\leq n-1} a^{ii}a^{jj}a_{ij,n}^2 +4D\sum_\i a_{ii,n}.
\end{split}
\end{align}
Then at $x_o$, by \eqref{MGrad}, we have
$$M_1+M_2 = N_1+N_2+N_3+N_4+N_5,$$
where
\begin{align*}
N_1=&-C \bigg(\sum_{i=2}^{n-1}a^{ii}a_{ii,n} \bigg)^2
-C\sum_{i=2}^{n-1}(a^{ii}a_{ii,n})^2
+4\sum_{i=2}^{n-1}\Big[D(a_{ii}-a_{11})-C\rho'u_nu_{nn}\Big]a^{ii}a_{ii,n},\\
N_2=&-A\sum_{2\leq i \leq n-1}(1+2a_{ii}a^{11})\cdot(a^{ii}a_{ii,
1})^2 -A\bigg(\sum_{2\leq i\leq n-1}a^{ii}a_{ii,1}\bigg)^2\\
&+4u_{n1}\sum_{2\leq i\leq n-1}
\Big[B(a_{ii}a^{11}-1)-A\rho'u_n\Big]\cdot(a^{ii}a_{ii,1}),\\
N_3=&\sum_{2\leq j\leq n-1}\bigg\{ -A(1+2a_{11}a^{jj})\cdot \bigg(
\sum_{2\leq i \leq n-1}a^{ii}a_{ii,j}\bigg)^2 -A\sum_{ \substack{2\leq i\leq n-1\\ i\neq j}}2a_{ii}a^{jj}\cdot(a^{ii}a_{ii, j})^2\\
&+4u_{nj}\sum_{2\leq i \leq n-1} \Big[Ba_{ii}a^{jj}-A\rho'u_n -(B+2A\rho'u_n)a_{11}a^{jj}\Big] \cdot(a^{ii}a_{ii,j})\\
&-A\sum_{2\leq i\leq n-1}(a^{ii}a_{ii,j})^2 \bigg\},
\end{align*}
and
\begin{align*}
N_4=&~2A\varphi_1\sum_\ii a^{ii}a_{ii,1}+4A\rho' u_n u_{n1} \varphi_1 +2A\sum_{2\leq i,j\leq n-1}(1+2a_{11}a^{jj}) a^{ii}a_{ii,j}\varphi_j\\
&+4A\rho' u_n\sum_{2\leq j\leq n-1}(1+2a_{11}a^{jj})u_{nj}\varphi_j
+4Ba_{11}\sum_{2\leq j\leq n-1} a^{jj}u_{nj}\varphi_j -A\varphi_1^2 -C\varphi_n^2\\
&+2\Big[C\big(\sum_{\ii}a^{ii}a_{ii,n}+2\rho' u_nu_{nn}\big) +2Da_{11}\Big]\varphi_n -A\sum_{2\leq j\leq n-1}(1+2a_{11}a^{jj}) \varphi_j^2,
\end{align*}

\begin{align*}
N_5=&~-u_n^2\sum_{\substack{1\leq i, j, k\leq n-1\\
i\neq j, j\neq k, k\neq i}}a^{ii}a^{jj}a_{ij, k}^2 -u_n^2\sum_{\substack {1\leq i,j\leq n-1\\i\neq j}} a^{ii} a^{jj} a_{ij,n}^2
-\Big(4A\rho'^2u_n^2+8B\rho'u_n\Big)u_{n1}^2\\
&-\sum_{2\leq j\leq n-1}\Big(4A\rho'^2u_n^2 +8A\rho'^2u_n^2a_{11}a^{jj}
+8B\rho'u_n a_{11}a^{jj}\Big)u_{nj}^2 -4C\rho'^2u_n^2u_{nn}^2\\
&-8D\rho' u_na_{11} u_{nn}.
\end{align*}
\end{lemma}
\begin{proof}
For the term $M_1$, we have
\begin{align}\label{W11-14}
\begin{split}
M_1 =&~M_{11}+M_{12}+M_{13}+M_{14},
\end{split}
\end{align}
where
\begin{align*}
M_{11}=&-A\sum_{\substack{1\leq i, j, k\leq n-1\\ i\neq j, j\neq k,
k\neq i}}a^{ii}a^{jj}a_{ij, k}^2,\quad
M_{12}=-A\sum_\i(a^{ii}a_{ii,i})^2,\\
M_{13}=&-2A\sum_{\substack{1\leq i, j \leq n-1\\ i\neq j}}
a^{ii}a^{jj}a_{ij, j}^2-A\sum_{\substack{1\leq i, j \leq n-1\\ i\neq j}}(a^{ii}a_{ii, j})^2, \\
M_{14}=&~4B\sum_{1\leq i,j\leq n-1}a^{jj}u_{nj}a_{ii, j}.
\end{align*}
By \eqref{MGrad},
\begin{align}\label{3grad'}
a^{11}a_{11,\a}= \varphi_\a-\sum_{2\leq i\leq n-1}a^{ii}a_{ii,\a} -2\rho' u_n u_{n\a},
\end{align}
hence
\begin{align}\label{W12}
\begin{split}
M_{12}=&-A(a^{11}a_{11,1})^2 -A\sum_{2\leq i\leq n-1}(a^{ii}a_{ii,i})^2\\
=&-A\bigg(\sum_{2\leq i \leq n-1}a^{ii}a_{ii,1}\bigg)^2 -4A\rho'
u_nu_{n1}\sum_{2\leq i\leq n-1}a^{ii}a_{ii,1} -4A\rho'^2u_n^2u_{n1}^2\\
&-A\sum_{2\leq i\leq n-1}(a^{ii}a_{ii,i})^2+2A\varphi_1\sum_\ii a^{ii}a_{ii,1} +4A\rho' u_nu_{n1}\varphi_1 -A\varphi_1^2,
\end{split}
\end{align}
and
\begin{align}\label{W13}
\begin{split}
M_{13}=&-A\sum_{2\leq j \leq n-1}(1+2a_{11}a^{jj})\cdot(a^{11}a_{11, j} )^2
-A\sum_{\substack{2\leq i \leq n-1\\ \j\\ i\neq j}} (1+2a_{ii}a^{jj}) \cdot(a^{ii}a_{ii,j})^2\\
=&-A\sum_{2\leq j\leq n-1} (1+2a_{11}a^{jj}) \cdot\bigg(\sum_{2\leq
i \leq n-1} a^{ii}a_{ii,j}\bigg)^2 -4A\rho'^2u_n^2\sum_{2\leq j \leq n-1}(1+2a_{11}a^{jj}) u_{nj}^2 \\
&-A\sum_{2\leq i \leq n-1}(1+2a_{ii}a^{11})\cdot(a^{ii}a_{ii, 1})^2
-A\sum_{\substack{2\leq i,j \leq n-1\\ i\neq
j}}(1+2a_{ii}a^{jj})\cdot(a^{ii}a_{ii, j})^2\\
&-4A\rho' u_n\sum_{2\leq i,j\leq n-1}(1+2a_{11}a^{jj}) u_{nj}a^{ii}a_{ii,j} +2A\sum_{2\leq i,j\leq n-1}(1+2a_{11}a^{jj}) a^{ii}a_{ii,j}\varphi_j\\
&+4A\rho' u_n\sum_{2\leq j\leq n-1}(1+2a_{11}a^{jj}) u_{nj}\varphi_j
-A\sum_{2\leq j\leq n-1}(1+2a_{11}a^{jj}) \varphi_j^2.
\end{split}
\end{align}
Making use of \eqref{3grad'} again, we can obtain
\begin{align}\label{W14}
\begin{split}
M_{14}=&~4B\sum_\i a^{ii}u_{ni}a_{11,i}
+4B\sum_{\substack{\i\\ 2\leq j \leq n-1}}a^{ii}u_{ni}a_{jj, i}\\
=&~4Bu_{n1}\sum_\ii (a^{11}-a^{ii})a_{ii,1} -4B\sum_{ 2\leq i, j
\leq n-1}a_{11}a^{ii}a^{jj}u_{ni}a_{jj,i}\\& +4B\sum_{ 2\leq i,j
\leq n-1}a^{ii}u_{ni}a_{jj, i}-8B\rho'u_na_{11}\sum_\i
a^{ii}u_{ni}^2 +4Ba_{11}\sum_\i a^{ii}u_{ni} \varphi_i.
\end{split}
\end{align}
By \eqref{W11-14}, \eqref{W12}--\eqref{W14},
\begin{align}\label{W1}
\begin{split}
M_1=&~-A\sum_{2\leq i \leq n-1}(1+2a_{ii}a^{11})\cdot(a^{ii}a_{ii,
1})^2 -A\bigg(\sum_{2\leq i\leq n-1}a^{ii}a_{ii,1}\bigg)^2\\
&+4u_{n1}\sum_{2\leq i\leq n-1}\Big[B(a_{ii}a^{11}-1)-A\rho'u_n\Big]\cdot(a^{ii}a_{ii,1})\\
&+\sum_{2\leq j\leq n-1}\bigg\{ -A(1+2a_{11}a^{jj})\cdot \bigg(
\sum_{2\leq i \leq n-1}a^{ii}a_{ii,j}\bigg)^2 -A\sum_{ \substack{2\leq i\leq n-1\\ i\neq j}}2a_{ii}a^{jj}\cdot(a^{ii}a_{ii, j})^2\\
&-A\sum_{2\leq i\leq n-1}(a^{ii}a_{ii,j})^2 +4u_{nj}\sum_{2\leq i \leq n-1}\big[Ba_{ii}a^{jj}-A\rho'u_n-(B+2A\rho'u_n)a_{11}a^{jj}\big]\cdot(a^{ii}a_{ii,j})\bigg\}\\
&-A\sum_{\substack{1\leq i, j, k\leq n-1\\ i\neq j, j\neq k, k\neq
i}} a^{ii}a^{jj}a_{ij, k}^2 -\sum_{2\leq j\leq n-1} \Big(4A\rho'^2
u_n^2+8A\rho'^2u_n^2a_{11}a^{jj} +8B\rho'u_n a_{11}a^{jj} \Big)u_{nj}^2\\
&-\Big(4A\rho'^2u_n^2+8B\rho'u_n\Big)u_{n1}^2+M_1(\n \varphi),
\end{split}
\end{align}
where
\begin{align*}
M_1(\n \varphi)=&~2A\varphi_1\sum_\ii a^{ii}a_{ii,1}+4A\rho' u_n
u_{n1} \varphi_1+2A\sum_{2\leq i,j\leq n-1}
(1+2a_{11}a^{jj}) a^{ii}a_{ii,j}\varphi_j\\
& +4A\rho' u_n\sum_{2\leq j\leq n-1}(1+2a_{11}a^{jj})u_{nj}\varphi_j
+4Ba_{11}\sum_{1\leq j\leq n-1} a^{jj}u_{nj}\varphi_j -A\varphi_1^2\\
&-A\sum_{2\leq j\leq n-1}(1+2a_{11}a^{jj}) \varphi_j^2.
\end{align*}

For the term $M_2$, we write it as
\begin{align}\label{W212}
M_2=M_{21}+M_{22},
\end{align}
where
\begin{align*}
M_{21}=-C\sum_{1\leq i,j\leq n-1}a^{ii}a^{jj}a_{ij,n}^2, \qquad
M_{22}=4D\sum_\i a_{ii,n}.
\end{align*}
Analogy, by \eqref{3grad'} we have
\begin{align}\label{W21}
\begin{split}
M_{21}=&-C\sum_{\substack{1\leq i,j\leq n-1\\i\neq j}
}a^{ii}a^{jj}a_{ij,n}^2
-C(a^{11}a_{11,n})^2 -C\sum_{2\leq i \leq n-1} (a^{ii}a_{ii,n})^2\\
=&-C\sum_{\substack{1\leq i,j\leq n-1\\i\neq j}} a^{ii}a^{jj} a_{ij,
n}^2 -4C\rho' u_nu_{nn}\sum_{_2\leq i\leq n-1}a^{ii}a_{ii,n}
-C\sum_{2\leq i \leq n-1} (a^{ii}a_{ii,n})^2\\
&-4C\rho'^2u_n^2u_{nn}^2 -C \bigg(\sum_{2\leq i \leq
n-1}a^{ii}a_{ii,n} \bigg)^2+2C\Big(\sum_{\ii}a^{ii}a_{ii,n}+2\rho'
u_nu_{nn}\Big)\varphi_n-C\varphi_n^2,
\end{split}
\end{align}
also
\begin{align}\label{W22}
\begin{split}
M_{22}=&~4Da_{11,n} +4D\sum_{2\leq i \leq n-1}a_{ii,n}\\
=&~-4D\sum_{2\leq i \leq n-1}a_{11}a^{ii}a_{ii,n} -8D\rho' u_na_{11}
u_{nn} +4D\sum_{2\leq i \leq n-1}a_{ii,n}+4Da_{11}\varphi_n.
\end{split}
\end{align}
Combining \eqref{W212}--\eqref{W22}, we obtain
\begin{align*}
M_2=&-C\sum_{\substack{1\leq i,j\leq n-1\\i\neq j}} a^{ii}a^{jj}
a_{ij, n}^2-C \bigg(\sum_{2\leq i \leq n-1}a^{ii}a_{ii,n} \bigg)^2
-C\sum_{2\leq i \leq n-1} (a^{ii}a_{ii,n})^2\\
& -4(C\rho' u_nu_{nn}+Da_{11})\sum_{2\leq i\leq n-1}a^{ii}a_{ii,n}
+4D\sum_{2\leq i \leq n-1}a_{ii,n}-4C\rho'^2u_n^2u_{nn}^2\\
&-8D\rho' u_na_{11}u_{nn} -C\varphi_n^2 +2\Big(C\sum_{\ii}a^{ii}a_{ii,n} +2C\rho' u_nu_{nn} +2Da_{11}\Big)\varphi_n.
\end{align*}

After the computation above, we denote by $N_1$ the terms involving
$a_{ii,n}(2\leq i \leq n-1)$; $N_2$ the terms involving $a_{ii,1}(2\leq i \leq
n-1)$; $N_3$ the terms involving $a_{ii,j}(2\leq i, j \leq n-1)$; $N_4$ the
terms involving $\n \varphi$ and $N_5$ the other terms. Then we complete the
proof.
\end{proof}

Now we state the following elementary calculus lemma, which had appeared in \cite{MOZ09}.
\begin{lemma}\label{maximum}
Let $\lambda\geq 0$, $\mu\in\mathbb{R}$, $b_i>0$ and $c_i\in\mathbb{R}$ for $2\leq i\leq n-1$. Define the quadratic polynomial
\begin{align*}
\mathcal {Q}(X_2, \cdots, X_{n-1})=-\sum_{2\leq i\leq n-1}b_iX_i^2
-\lambda\bigg(\sum_{2\leq i\leq n-1}X_i\bigg)^2 +4\mu\sum_{2\leq
i\leq n-1}c_iX_i.
\end{align*}
Then we have $$\mathcal {Q}(X_2,\cdots, X_{n-1})\leq 4\mu^2\Gamma,$$
where
$$\Gamma=\sum_{2\leq i\leq n-1}\frac{c_i^2}{b_i}-\lambda\bigg(1+\lambda
\sum_{2\leq i\leq n-1} \frac{1}{b_i}\bigg)^{-1}\bigg(\sum_{2\leq
i\leq n-1}\frac{c_i}{b_i}\bigg)^2.$$
\end{lemma}

\section{Proof of Theorem \ref{MinimalThm}.}

In this section, we shall study the following equation
\begin{align}\label{Mean1}
{\rm div}(\frac{\n u}{\sqrt{1+|\n u|^2}})=0\qquad {\text {in}}\ \Omega \subset \mathbb{R}^n,
\end{align}
and prove the Theorem \ref{MinimalThm}.

Denote
\begin{align}\label{MFab1}
F^{\alpha \beta}=(1+|\nabla u|^2)\delta_{\alpha \beta}-u_{\alpha}u_{\beta},
\end{align}
then equation \eqref{Mean1} is reduced  to
\begin{equation}\label{Mean2}
F^{\alpha \beta}u_{\alpha \beta}=0.
\end{equation}

As we have mentioned in the last section, set
$$\varphi=\rho(|\nabla u|^2)+ \log K(x).$$
For suitable choice of $\rho$, we will derive the following elliptic differential inequality
\begin{align*}
\sum_{1\leq \alpha, \beta \leq n}F^{\alpha \beta}\varphi_{\alpha \beta} \leq 0
\qquad \mod \nabla \varphi \quad \text{in} \quad \Omega,
\end{align*}
where we modify the terms of $\nabla \varphi$ with locally bounded
coefficients. We shall complete the calculation at the fixed point $x_o$. As in the last section we may assume that $u_i(x_o)=0 (1\leq i\leq n-1)$ and
$u_n(x_o) =|\nabla u|>0$. And we can further assume that the matrix
$(u_{ij}(x_o)) (1\leq i, j\leq n-1)$ is diagonal and $u_{ii}(x_o) < 0$.

In this section, all the calculation will be worked at $x_o$.

By \eqref{MFab1}, we have
\begin{equation}\label{Fab}
\begin{split}
F^{ii}=1+u_n^2~( 1 \leq i \leq n-1), \qquad F^{nn}=1,\qquad
F^{\alpha\beta}=0 \ \ \text{for} \ \ \alpha \neq \beta,
\end{split}
\end{equation}
and
\begin{equation}\label{Fab1}
\begin{split}
F^{\a\b}_{,\gamma}&=2u_nu_{n\gamma}\delta_{\alpha \beta}-u_{\a \gamma} u_{\b} -u_{\a}u_{\b \gamma},\\
F^{\a\b}_{,ii}&=2\sum_{1\leq \gamma \leq n}u_{\gamma i}^2 \delta_{\alpha \beta} +2u_nu_{iin}\delta_{\alpha \beta}-u_{\a ii} u_{\b} -u_{\a }u_{\b ii} -2u_{\a i} u_{\b i}.
\end{split}
\end{equation}

It is also easy to check
\begin{align}\label{Munn}
u_{nn}=u_n(1+u_n^2)\sigma_1,\qquad \Delta u =  u_n^3\sigma_1,\quad \text{where} \quad \sigma_1=\sum_\i a_{ii}.
\end{align}

Recalling the following formula we obtained in \eqref{H1234},
\begin{align}\label{H1234a}
\sum_{1\leq \alpha, \beta \leq n} F^{\alpha\beta} \varphi_{\alpha \beta} =L_1+L_2+L_3+L_4,
\end{align}
we will treat the terms $L_1$, $L_2$, $L_3$ and $L_4$ above respectively.

Let us deal with the term $L_1$ in \eqref{H1234a} at first. By the equation \eqref{Mean2}, we have
\begin{align*}
L_1=~L_{11} + L_{12}+L_{13},
\end{align*}
where
\begin{align*}
L_{11}=&~6u_n^{-2}\sum_{\i} a^{ii}u_{ni} \sum_{1\leq\a,\b \leq n} F^{\a\b} u_{\a\b i},\\
L_{12}=&~[2\rho 'u_n-(n-1)u_n^{-1}] \sum_{1\leq \a, \b \leq n} F^{\a\b} u_{\a\b n}, \\
L_{13}=&-{u_n}^{-1}\sum_{\i}a^{ii}\sum_{1\leq \alpha,\beta\leq n}F^ {\alpha \beta} u_{ \alpha \beta ii}.
\end{align*}
Then by \eqref{Fab1}, we have
\begin{align}\label{MH12}
\begin{split}
L_{11} =&-6u_n^{-2}\sum_{\i}a^{ii}u_{ni}\sum_{1\leq \a, \b \leq
n}F^{\a\b}_{i}u_{\a\b} \\
=&~6u_n^{-2}\sum_{\i}a^{ii}u_{ni}\bigg[2u_n
u_{ni}(u_n\sigma_1  + u_{ii})\bigg]\\
=&~12\sigma_1\sum_{\i}a^{ii}u_{ni}^2-12\sum_{\i}u_{ni}^2,
\end{split}
\end{align}
and
\begin{align}\label{MH14}
\begin{split}
L_{12} =&-[2\rho 'u_n-(n-1)u_n^{-1}] \sum_{1\leq \a, \b \leq n} F^{\a\b}_n u_{\a\b}\\
=&~[2\rho 'u_n-(n-1)u_n^{-1}] \bigg(2u_n^2\sigma_1u_{nn} +2u_n\sum_\j u_{nj}^2 \bigg).
\end{split}
\end{align}

Now we shall calculate the term $L_{13}$. By differentiating \eqref{Mean2} twice with respect to $x_i$, we have
\begin{equation*}
\sum_{1\leq \a, \b \leq n} F^{\a\b}u_{\a\b ii} = -\sum_{1 \leq \a,\b \leq n} F^{\a\b}_{,ii}u_{\a\b} - 2 \sum_{1 \leq \a,\b \leq n} F^{\a\b}_{,i}u_{\a\b i}.
\end{equation*}
By \eqref{Fab1},
\begin{align*}
-u_n^{-1}\sum_{1\leq\a, \b\leq n} F^{\a\b}u_{ii\a\b}
=&~u_n^{-1}\sum_{1 \leq \a,\b \leq n}F^{\a\b}_{,ii}u_{\a\b} +2u_n^{-1}
\sum_{1 \leq \a,\b \leq n}F^{\a\b}_{,i}u_{\a\b i}\\
=&\ 2u_n^{-1}\Delta u \sum_{1\leq \a \leq n}u_{i\a}^2 +2\Delta u u_{nii} -2\sum_{1\leq \a \leq n}u_{\a ii}u_{n\a}\\
&-2u_n^{-1}\sum_{1\leq \a\leq n}u_{\a i} u_{\b i}u_{\a\b} +4
u_{ni}(\Delta u)_i -4\sum_{1\leq \a \leq n}u_{\a i}u_{ni\a}\\
=&L_{13a}+L_{13b}+L_{13c},
\end{align*}
where
\begin{align*}
L_{13a}=&2u_n^{-1}\Delta u \sum_{1\leq \a \leq n}u_{i\a}^2-2u_n^{-1} \sum_{1\leq \a\leq n}u_{\a i} u_{\b i}u_{\a\b} ,\quad L_{13b} =2\Delta u u_{nii} -2\sum_{1\leq \a \leq n}u_{\a ii}u_{n\a},\\
L_{13c}=&4 u_{ni}(\Delta u)_i -4\sum_{1\leq \a \leq n}u_{\a i}u_{ni\a}.
\end{align*}

Hence, for the minimal surface equation, we have
\begin{align}\label{MH111ok}
\begin{split}
\sum_{\i}a^{ii}L_{13a}=-2\sigma_1\sum_{\i}a^{ii}u_{ni}^2+2u_n\Delta
u\sigma_1+4\sum_{\i}u_{ni}^2+2u_n^{2}\sum_{\i}a_{ii}^2.
\end{split}
\end{align}

Thanks to \eqref{uiia} and \eqref{3grad'}, we have
\begin{align}\label{MH112ok}
\begin{split}
\sum_{\i}a^{ii}L_{13b}=&~2\sum_{\i}a^{ii}u_{iin}\sum_{\j}u_{jj}-2 \sum_{1\le i,j \le n-1} a^{ii}u_{iij}u_{nj}\\
=&-2u_n\sigma_1\sum_{\i}a^{ii}(-u_na_{ii,n}+2u_n^{-1}u_{ni}^2 +u_n^{-1} u_{nn} u_{ii})\\
&-2 \sum_{1\le i,j\le n-1}u_{nj}a^{ii}(-u_na_{ii,j}+2u_n^{-1}u_{ni}u_{ij} +u_n^{-1}u_{nj}u_{ii})\\
=&~[2(n-1)-4\rho'u_n^2]u_nu_{nn}\sigma_1+[2(n+1)-4\rho'u_n^2]\sum_{\i}u_{ni}^2\\
&-4\sigma_1\sum_{\i}a^{ii}u_{ni}^2+2u_n^{2}\sigma_1\varphi_n+2u_n\sum_{\j}u_{nj}\varphi_j,
\end{split}
\end{align}
and
\begin{align}\label{MH113ok}
\begin{split}
\sum_{\i}a^{ii}L_{13c}=&~4\sum_{\i}a^{ii}u_{ni}\sum_{\j}u_{jji}+4u_n\sum_{\i}u_{iin}\\
=&~4\sum_{\i}a^{ii}u_{ni}\sum_{\j}(-u_na_{jj,i}+2u_n^{-1}u_{nj}u_{ji} +u_n^{-1} u_{ni} u_{jj})\\
&+4u_n\sum_{\i}(-u_na_{ii,n}+2u_n^{-1}u_{ni}^2 +u_n^{-1}u_{nn}u_{ii})\\
=&-4\sigma_1\sum_{\i}a^{ii}u_{ni}^2 -4u_n^2\sum_{\i}a_{ii,n}-4u_nu_{nn} \sigma_1\\
&-4u_n\sum_{1\le i,j \le n-1}a^{ii}u_{ni}a_{jj,i}.
\end{split}
\end{align}

By \eqref{MH111ok}--\eqref{MH113ok}, it yields
\begin{align}\label{MH11ok}
\begin{split}
L_{13}=&-u_n^{-1}\sum_{\i}a^{ii}\sum_{1\leq \alpha,\beta\leq
n}F^ {\alpha\beta}u_{ii, \alpha \beta}\\
=&-4u_n^2\sum_{\i}a_{ii,n}-4u_n\sum_{1\le i,j \le
n-1}a^{ii}u_{ni}a_{jj,i}+2u_n^{2}\sum_{\i}a_{ii}^2\\
&+\big(2n-6-4\rho'u_n^2\big)u_nu_{nn}\sigma_1+2u_n\Delta
u \sigma_1+\big(2n+6-4\rho'u_n^2\big)\sum_{\i}u_{ni}^2\\
&-10\sigma_1\sum_{\i}a^{ii}u_{ni}^2 +2u_n\sum_{\j}u_{nj}\varphi_j
+2u_n^{2}\sigma_1\varphi_n.
\end{split}
\end{align}

Combining \eqref{MH12}, \eqref{MH14} and \eqref{MH11ok}, we obtain
\begin{align}\label{MH1ok}
\begin{split}
L_{1}=& -4u_n^2\sum_{\i}a_{ii,n}-4u_n\sum_{1\le i,j \le
n-1}a^{ii}u_{ni}a_{jj,i}-4u_nu_{nn}\sigma_1\\&+2u_n\Delta u \sigma_1
+2\sigma_1\sum_{\i}a^{ii}u_{ni}^2 +2u_n^{2}\sum_{\i}a_{ii}^2\\& -4\sum_{\i}u_{ni}^2+2u_n^{2}\sigma_1\varphi_n +2u_n\sum_{\j}u_{nj}\varphi_j.
\end{split}
\end{align}

For the term $L_2$ in \eqref{H1234a}, applying \eqref{uiia} and \eqref{Fab}, we have
\begin{align}\label{MH2}
\begin{split}
L_2 =&- (1+u_n^2)\sum_{1\leq i,j,k\leq n-1} a^{ii}a^{jj}a_{ij, k}^2-
\sum_{1\leq i,j\leq n-1} a^{ii}a^{jj}a_{ij, n}^2\\
&-4u_n^{-1}(1+u_n^2)\sum_{\i}u_{iin} -4(1+u_n^{-2})\sum_{1\leq
i,j\leq n-1}a^{ii}u_{ni}u_{jj i}\\
=&- (1+u_n^2)\sum_{1\leq i,j,k\leq n-1} a^{ii}a^{jj}a_{ij, k}^2-
\sum_{1\leq i,j\leq n-1} a^{ii}a^{jj}a_{ij,n}^2\\
&+4(1+u_n^2)\sum_{\i}a_{ii,n}+4u_n^{-1}(1+u_n^{2})\sigma_1u_{nn}\\
&+4u_n^{-1}(1+u_n^2)\sum_{1\le i,j \le n-1}a^{ii}u_{ni}a_{jj,i}
+4(1+u_n^{-2})\sigma_1\sum_{\i}a^{ii}u_{ni}^2.
\end{split}
\end{align}

For the term $L_3$ and $L_4$, in a similar way, we can obtain
\begin{align}\label{MH34}
\begin{split}
L_3=&~\big[2\rho 'u_n^{2}-(n+1)\big](1+u_n^{2})\sum_{\j}a_{jj}^2+[4\rho''
u_n^{4}+6\rho'u_n^2](1+u_n^{2})\sigma_1^2\\
&-6u_n^{-2}(1+u_n^{2})\sigma_1\sum_{\i}a^{ii}u_{ni}^2-2u_n^{-1}(1+u_n^{2})u_{nn}\sigma_1\\
&+\bigg[(6\rho'u_n^2 + 4\rho''u_n^4+4)+ (8\rho'u_n^2 + 4\rho''u_n^4
-n+3)u_n^{-2}\bigg]\sum_{\i}u_{ni}^2,\\
L_4=&\,-2u_n^{-1}u_{nn}\varphi_{n}-2u_n^{-1}(1+u_n^{2})\sum_{\i}u_{ni}\varphi_{i}.
\end{split}
\end{align}

With \eqref{Munn} in hand, combining \eqref{MH1ok}---\eqref{MH34}, we finally get
\begin{align}\label{Mfinish}
\begin{split}
\sum_{1\leq\a,\b\leq n}F^{\a\b}\varphi_{\a\b}
=&-(1+u_n^2)\sum_{1\leq i,j,k\leq n-1} a^{ii}a^{jj}a_{ij,k}^2
-\sum_{1\leq i,j \leq n-1} a^{ii}a^{jj}a_{ij,n}^2 +4\sum_{\i}a_{ii,n}\\
&+4u_n^{-1}\sum_{1\le i,j\le n-1}a^{ii}u_{ni}a_{jj,i}
+\bigg[(2\rho'u_n^2-n-1)+(2\rho'u_n^2-n+1)u_n^2\bigg]\sum_{\i}a_{ii}^2\\
& +\bigg[(6\rho'u_n^2 + 4\rho''u_n^4)+ (8\rho'u_n^2
+ 4\rho''u_n^4 -n+3)u_n^{-2}\bigg]\sum_\i u_{ni}^2\\
&+\bigg[4\rho''u_n^4(1+u_n^2)^2+6\rho'u_n^2(1+u_n^2)^2
+2\bigg]\sigma_1^2-2u_n^{-2}\sigma_1\sum_\i a^{ii}u_{ni}^2\\
&-2u_n^{-1}\sum_\j u_{nj}\varphi_j-2\sigma_1\varphi_n.
\end{split}
\end{align}

Now, we shall try a test function to estimate the Gaussian curvature of the level sets of the minimal graph. By setting $\rho(t)= \t[\log t -\log(1+t)]$ in \eqref{Mfinish}, we obtain
\begin{align}\label{Mp1}
\begin{split}
\sum_{1\leq\a,\b\leq n}F^{\a\b}\varphi_{\a\b}
=&-(1+u_n^2)\sum_{1\leq i,j,k\leq n-1} a^{ii}a^{jj}a_{ij,k}^2
-\sum_{1\leq i,j \leq n-1} a^{ii}a^{jj}a_{ij,n}^2 +4\sum_{\i}a_{ii,n}\\
&+4u_n^{-1}\sum_{1\le i,j\le n-1}a^{ii}u_{ni}a_{jj,i}
+\bigg[(2\t-n-1)-(n-1)u_n^2\bigg]\sum_{\i}a_{ii}^2\\
& +\bigg[2\t(1-u_n^2)(1+u_n^2)^{-2}+ \Big(4\t(1+u_n^2)^{-2} -n+3\Big)u_n^{-2}\bigg]\sum_\i u_{ni}^2\\
&+\bigg[2\t(1-u_n^2)+2\bigg]\sigma_1^2-2u_n^{-2}\sigma_1 \sum_\i
a^{ii}u_{ni}^2\\& -2u_n^{-1}\sum_\j u_{nj}\varphi_j
-2\sigma_1\varphi_n.
\end{split}
\end{align}

Let us solve the 2-dimensional case at first. Now the formula \eqref{Mp1} reduces to
\begin{align}\label{Mpn2}
\begin{split}
\sum_{1\leq\a,\b\leq 2}F^{\a\b}\varphi_{\a\b} =&-(1+u_2^2)(a^{11}a_{11,1})^2 -(a^{11}a_{11,2})^2 +4a_{11,2}+4u_2^{-1} u_{21}a^{11}a_{11,1}\\
&+\bigg[2\t(1-u_2^2)(1+u_2^2)^{-2}+ \Big(4\t(1+u_2^2)^{-2} -1\Big)u_2^{-2}\bigg] u_{21}^2\\
&+\bigg[4\t-1-(2\t+1)u_2^2\bigg]a_{11}^2 -2u_2^{-1}u_{21}\varphi_1 -2\sigma_1 \varphi_2.
\end{split}
\end{align}
By \eqref{MGrad}, we have
$$a^{11}a_{11,\alpha}=\varphi_\alpha-2\t u_2^{-1}(1+u_2^2)^{-1}u_{2\alpha}, \quad{\text{for}}\ \alpha=1,2.$$
For $\t=-\frac 1 2$, it reads
\begin{align}\label{Mpn2ok}
\begin{split}
\sum_{1\leq\a,\b\leq 2}F^{\a\b}\varphi_{\a\b} =-(1+u_2^2)\varphi_1^2
-\varphi_2^2,
\end{split}
\end{align}
which shows the validity of the Theorem \ref{MinimalThm} in 2-dimensional case via the strong minimum (maximum) principle.

In the following, we come to deal with the case $n\geq3$. Let
\begin{align*}
\sum_{1\leq\a,\b\leq n}F^{\a\b}\varphi_{\a\b} = P_1+P_2 +P_3,
\end{align*}
where
\begin{align*}
P_1=&-(1+u_n^2)\sum_{1\leq i,j,k\leq n-1} a^{ii}a^{jj}a_{ij,k}^2 +4u_n^{-1}\sum_{1\le i,j\le n-1}a^{ii}u_{ni}a_{jj,i}\\
&-\sum_{1\leq i,j \leq n-1} a^{ii}a^{jj}a_{ij,n}^2 +4\sum_{\i}a_{ii,n},
\end{align*}
and
\begin{align*}
P_2=&~\bigg[(2\t-n-1)-(n-1)u_n^2\bigg]\sum_{\i}a_{ii}^2+\bigg[2\t(1-u_n^2)+2\bigg]\sigma_1^2\\
&+\bigg[2\t(1-u_n^2)(1+u_n^2)^{-2}+ \Big(4\t(1+u_n^2)^{-2} -n+3\Big)u_n^{-2} \bigg] \sum_\i u_{ni}^2\\
&-2u_n^{-2}\sigma_1 \sum_\i a^{ii}u_{ni}^2,\\
P_3=& -2u_n^{-1}\sum_\j u_{nj}\varphi_j-2\sigma_1\varphi_n.
\end{align*}

To deal with the term $P_1$, we can set $A=1+u_n^2,\ B=u_n^{-1},\ C=D=1$ in Lemma \ref{Form2}. Let us denote by $(Q_1)$ the terms involving $a_{ii,n}(2\leq i \leq n-1)$; $(Q_2)$ the terms involving $a_{ii,1}(2\leq i \leq n-1)$; $(Q_3)$ the terms involving $a_{ii,j}(2\leq i, j \leq n-1)$; $(Q_4)$ the terms involving $\n \varphi$ and $(Q_5)$ all of the rest terms. More precisely, we have
\begin{align}\label{biaodashi'}
\sum_{1\leq\a,\b\leq n}F^{\a\b}\varphi_{\a\b} = Q_1+Q_2+Q_3+Q_4+Q_5,
\end{align}
where
\begin{align*}
Q_1=&- \bigg(\sum_{2\leq i \leq n-1}a^{ii}a_{ii,n} \bigg)^2 -\sum_{2\leq i \leq n-1} (a^{ii}a_{ii,n})^2+4\sum_{_2\leq i\leq n-1} (a_{ii}-a_{11} -\t\sigma_1) a^{ii}a_{ii,n},\\
Q_2=&-(1+u_n^2)\sum_{2\leq i \leq n-1} (1+2a_{ii}a^{11})\cdot(a^{ii}a_{ii, 1})^2 -(1+u_n^2)\bigg( \sum_{2\leq i\leq n-1}a^{ii}a_{ii,1}\bigg)^2\\
&+4u_n^{-1}u_{n1}\sum_{2\leq i\leq n-1}\Big[(a_{ii}a^{11}-1)-\t\Big] \cdot(a^{ii}a_{ii,1}),\\
Q_3=&\sum_{2\leq j\leq n-1}\bigg\{ -(1+u_n^2)(1+2a_{11}a^{jj}) \cdot \bigg( \sum_{2\leq i \leq n-1}a^{ii}a_{ii,j}\bigg)^2\\
&+4u_n^{-1}u_{nj}\sum_{2\leq i \leq n-1}\Big[ a_{ii}a^{jj}-\t-(1+2\t)a_{11} a^{jj}\Big] \cdot(a^{ii}a_{ii,j})\\
\end{align*}
\begin{align*}
&-(1+u_n^2)\sum_{\substack{2\leq i\leq n-1\\ i\neq j}}2a_{ii}a^{jj}
\cdot(a^{ii}a_{ii, j})^2 -(1+u_n^2)\sum_{2\leq i\leq n-1} (a^{ii}
a_{ii,j})^2 \bigg\},\\
Q_4=&~2(1+u_n^2)\varphi_1\sum_\ii a^{ii}a_{ii,1} +4\t u_n^{-1} u_{n1} \varphi_1+2(1+u_n^2)\sum_{2\leq i,j\leq n-1} (1+2a_{11}a^{jj}) a^{ii}a_{ii,j} \varphi_j\\
&+4\t  u_n^{-1}\sum_{2\leq j\leq n-1}(1+2a_{11}a^{jj})u_{nj}\varphi_j +4u_n^{-1}a_{11}\sum_{2\leq j\leq n-1} a^{jj}u_{nj}\varphi_j -(1+u_n^2) \varphi_1^2 -\varphi_n^2\\
&+2\Big[\big(\sum_{\ii}a^{ii}a_{ii,n}+2\t(1+u_n^2)^{-1} u_n^{-1}u_{nn}\big) +2a_{11}\Big]\varphi_n\\
&-(1+u_n^2)\sum_{2\leq j\leq n-1}(1+2a_{11}a^{jj}) \varphi_j^2+P_3,
\end{align*}
and
\begin{align*}
Q_5=&~P_2-(1+u_n^2)\sum_{\substack{1\leq i, j, k\leq n-1\\ i\neq j, j\neq k, k\neq i}}a^{ii}a^{jj}a_{ij, k}^2 -\sum_{\substack{1\leq i,j\leq n-1\\ i\neq j}} a^{ii} a^{jj} a_{ij,n}^2\\
&-u_n^{-2}(1+u_n^2)^{-1}\sum_{2\leq j\leq n-1}\Big(4\t^2+8\t^2a_{11}a^{jj} +8\t a_{11}a^{jj}\Big)u_{nj}^2\\
&-(4\t^2+8\t)u_n^{-2}(1+u_n^2)^{-1}u_{n1}^2 -4\t^2\sigma_1^2-8\t a_{11} \sigma_1.
\end{align*}

By Lemma \ref{maximum}, we will maximize the terms $Q_1$, $Q_2$ and $Q_3$ for appropriate parameters. At first, let us examine the term $Q_1$. For $\ii$, set $X_i=a^{ii}a_{ii,n}$, $\lambda=1$, $\mu=1$, $b_i=1$ and $c_i =a_{ii} -a_{11}-\t \sigma_1$. By Lemma \ref{maximum}, we have
\begin{align}\label{Tunlast}
\begin{split}
Q_1\leq&~4\bigg[\sum_\ii c_i^2-\frac{1}{n-1}\bigg(\sum_\ii c_i\bigg)^2\bigg]\\
=&~4\sum_\j a_{jj}^2 +\frac{4(n-2)}{(n-1)} \t^2\sigma_1^2 -\frac{8}{n-1}\t \sigma_1^2 -\frac{4}{n-1} \sigma_1^2+8\t a_{11}\sigma_1.
\end{split}
\end{align}

For the term $(1+u_n^2)^{-1}Q_2$, in Lemma \ref{maximum}, set $X_i=a^{ii}a_{ii,1}$, $\lambda=1$, $\mu=u_{n1}u_n^{-1}(1+u_n^2)^{-1}$, $b_i=1+2a_{ii}a^{11}$ and $c_i=a_{ii}a^{11}-1-\t$ for $\ii$. Also by Lemma \ref{maximum}, we have
\begin{align*}
\frac{Q_2}{(1+u_n^2)}\leq&~\frac{4u_{n1}^2}{u_n^{2}(1+u_n^2)^{2}}\Gamma_1\\
=&~\frac{4u_{n1}^2}{u_n^{2}(1+u_n^2)^{2}}\Bigg[\sum_\ii\frac{c_i^2}{b_i}
-\bigg(1+\sum_\ii\frac{1}{b_i}\bigg)^{-1} \bigg(\sum_\ii \frac{c_i}{b_i} \bigg)^2 \Bigg].
\end{align*}
Let us simplify $\Gamma_1$. By denoting $$\beta_i=\frac{1}{b_i},$$ we have
$$a_{ii}a^{11}=\frac{1}{2\beta_i}-\frac12, \qquad c_i=\frac{1}{2\beta_i} -\frac32-\t.$$ Hence,
\begin{align*}
\Gamma_1=&\sum_\ii \beta_i\big(\frac{1}{2\beta_i}-\frac32-\t\big)^2
-\bigg(1+\sum_\ii \beta_i\bigg)^{-1}\bigg[\sum_\ii \beta_i\big( \frac{1}{2\beta_i}-\frac32-\t\big)\bigg]^2\\
=&~\frac14\sum_\ii \frac{1}{\beta_i} -\bigg(1+\sum_\ii \beta_i\bigg)^{-1} \big(\frac{n+1}{2}+\t\big)^2 +\big(\frac32+\t\big)^2,
\end{align*}
since
$$1 \le 1+\sum_\ii\beta_i \le n-1,$$
it follows that
\begin{align*}
\Gamma_1 \leq&~\frac14\sum_\ii \frac{1}{\beta_i} -\frac{1}{n-1}
\big(\frac{n+1}{2}+\t\big)^2 +\big(\frac32+\t\big)^2\\
=&~\frac{n-2}{n-1}(1+\t)^2+\frac14(2\sigma_1a^{11}-2).
\end{align*}
Then we have
\begin{align}\label{Tun1last}
Q_2\leq u_n^{-2}(1+u_n^2)^{-1} \bigg[\frac{4(n-2)}{n-1}(1+\t)^2 +2\sigma_1a^{11} -2\bigg]u_{n1}^2.
\end{align}

Let us estimate the term $(1+u_n^2)^{-1}Q_3$. We will apply Lemma \ref{maximum} for every $j$ fixed. For $\ii$, set $X_i=a^{ii}a_{ii,j}$, $\lambda=1 +2a_{11} a^{jj}$, $\mu=u_{nj}u_n^{-1}(1+u_n^2)^{-1}$, $b_i=1+2a_{ii}a^{jj}\ (i\neq j),
b_j=1$ and $c_i=a_{ii}a^{jj}-\t -(1+2\t)a_{11}a^{jj}$.

By Lemma \ref{maximum}, we have
$$(1+u_n^2)^{-1}Q_3\leq 4 u_n^{-2}(1+u_n^2)^{-2}\sum_{2\leq j\leq n-1} \Gamma_ju_{nj}^2,$$
where
$$\Gamma_j=c_j^2+\sum_{\substack{\ii\\i\neq j}}\frac{c_i^2}{b_i} -\bigg( \frac1\lambda+1+\sum_{\substack{\ii\\i\neq j}}\frac{1}{b_i} \bigg)^{-1} \bigg(c_j+\sum_{\substack{\ii\\i\neq j}}\frac{c_i}{b_i}\bigg)^2.$$
Also denoting
$$\beta_i=\frac{1}{b_i}\ (i\neq j),$$
we have
$$a_{ii}a^{jj}=\frac{1}{2\beta_i}-\frac12,\qquad c_i=\frac{1}{2\beta_i} -\delta,\ \  \mbox{ where}\ \  \delta =\frac12+\t +(1+2\t)a_{11}a^{jj}.$$
Noticed that
$$c_j=\frac32-\delta,\qquad \frac\delta\lambda=\frac12+\t,$$
we obtain
\begin{align*}
\Gamma_j=&~c_j^2+\sum_{\substack{\ii\\i\neq
j}}\beta_i(\frac{1}{2\beta_i}-\delta)^2-\bigg(\frac1\lambda+1
+\sum_{\substack{\ii\\i\neq j}}\beta_i\bigg)^{-1}\bigg[c_j
+\sum_{\substack{\ii\\i\neq j}}\beta_i(\frac{1}{2\beta_i}-\delta)
\bigg]^2\\
=&~\frac14\sum_{\substack{\ii\\i\neq j}}\frac{1}{\beta_i}
-\bigg(\frac1\lambda+1 +\sum_{\substack{\ii\\i\neq
j}}\beta_i\bigg)^{-1}(\frac
n2+\frac\delta\lambda)^2+\frac94+\frac{\delta^2}{\lambda}\\
=&~\frac14\sum_{\substack{\ii\\i\neq j}}
\frac1{\beta_i}-\bigg(\frac1\lambda+1+\sum_{\substack{\ii\\i\neq
j}}\beta_i\bigg) ^{-1}(\frac{n+1}{2}+\t)^2+\frac94 +(\frac12+\t)\delta.
\end{align*}
Obviously,
$$1 \leq \frac1\lambda+1+\sum_{\substack{\ii\\i\neq j}} \beta_i\leq n-1,$$
hence
\begin{align*}
\Gamma_j\leq&~\frac14\sum_{\substack{\ii\\i\neq j}} \frac1{\beta_i} -\frac{1}{n-1}(\frac{n+1}{2}+\t)^2+\frac94 +(\frac12+\t)\delta\\
=&~\frac{n-2}{n-1}\t^2-\frac{2}{n-1}\t +\frac{n-3}{2(n-1)} +\frac12\sigma_1a^{jj} +2\t^2a_{11}a^{jj} +2\t a_{11}a^{jj}.
\end{align*}
Therefore, we have
\begin{align}\label{Tunjlast}
Q_3\leq \frac {1}{u_n^{2}(1+u_n^2)}\sum_{2\leq j \leq n-1} \bigg[ \frac{4n-8}{n-1}\t^2 -\frac{8}{n-1}\t +\frac{2n-6}{n-1} +(2\sigma_1 +8\t^2a_{11} +8\t a_{11})a^{jj}\bigg]u_{nj}^2.
\end{align}

If we let
\begin{align*}
q_1(\t)=&~2\t-n+3+(-n+1)u_n^2,\\
q_2(\t)=&-\frac{4}{n-1}\t^2+(2-\frac{8}{n-1})\t+2-\frac{4}{n-1}-2\t u_n^2,\\
q_{3,j}(\t)=&\frac{1}{u_n^{2}(1+u_n^2)}\bigg[-\frac{4}{n-1}\t^2+(4-\frac{8}{n-1})\t
+5-n-\frac{4}{n-1}+(-2\t-n+3-2\sigma_1a^{jj})u_n^2\bigg],
\end{align*}
then collecting \eqref{biaodashi'}--\eqref{Tunjlast}, we finally obtain
\begin{align}\label{jixiaozuihou}
\begin{split}
\sum_{1\leq\a,\b\leq n}F^{\a\b}\varphi_{\a\b}\le ~ q_1(\t) \sum_{\j}a_{jj}^2+q_2(\t)\sigma_1^2+\sum_{\j}q_{3,j}(\t)u_{nj}^2
\quad \mod \n \varphi,
\end{split}
\end{align}
where we modify the terms of $\nabla \varphi$ with locally bounded coefficients. By a simple observation, one can see that a sufficient condition to validate
\begin{align*}
\sum\limits_{1\leq\a,\b\leq n}F^{\a\b}\varphi_{\a\b} \leq0\ \mod\nabla \varphi,
\end{align*}
is
\begin{align}\label{Condition}
 \left\{ \begin{array}{l}\vspace{0.2cm}
 q_1(\t) +q_2(\t) \leq0,\\ \vspace{0.2cm}\displaystyle
 q_1(\t) +(n-1)q_2(\t) \leq0,\\
 \displaystyle q_{3,j}(\t) \leq0.\\
 \end{array} \right.
\end{align}
Solving $q_1(\t)+q_2(\t)\leq0$, we have
\begin{equation*}
\t\geq\frac{1-n}{2}.
\end{equation*}
And $q_1(\t)+(n-1)q_2(\t)\leq 0$ implies
\begin{equation*}
\t=-\frac12 \quad \text{or}\quad \t\geq\frac{n-3}{2},
\end{equation*}
also from $q_{3,j}(\t)\leq0$, we get
\begin{equation*}
\t\geq \frac{1-n}{2}.
\end{equation*}
Therefore by solving the inequalities in \eqref{Condition}, we finally obtain $\t=-\dfrac 1 2$ or $\t \ge \dfrac{n-3}{2}$. So we complete our proof of Theorem \ref{MinimalThm}.\qed

In the 2-dimensional case, from the above proof we have the following observation.
\begin{Cor}\label{2-dim}
Let $\Omega\subset\mathbb{R}^2$ be a smooth bounded domain and $u\in C^4(\Omega) \cap C^{2}(\bar{\Omega})$ satisfies
\begin{align*}
{\rm div}\bigg(\frac{\n u}{\sqrt{1+|\n u|^2}}\bigg)=0\qquad {\text {in}}\ \Omega.
\end{align*}
Assume $|\nabla u|\neq 0$ in $\Omega$. Let $k$ be the curvature of the level curve, then the function $\bigg(\dfrac{|\n u|^2}{1+|\n u|^2}\bigg)^{-\frac12}k$ is a harmonic function with respect to the Laplace-Beltrami operator on the graph of $u$.
\end{Cor}
\begin{proof}
Let $\psi = \bigg(\dfrac{|\n u|^2}{1+|\n u|^2}\bigg)^{-\frac12}k$. Similar  as
the calculation leading to \eqref{Mpn2ok}, we have
\begin{align*}
\sum_{1\leq\a,\b\leq 2}F^{\a\b}\psi_{\a\b} =0 \qquad {\text {in}}\ \Omega.
\end{align*}
In other words, the function $\psi$ is a harmonic function  with respect to the Laplace-Beltrami operator on the graph of $u$(see \cite{GT98}).
\end{proof}
\begin{Rem}\label{2d}
Let $u$ be a 2-dimensional  harmonic function with no critical points in domain and $k$ be the curvature of the level curve of $u$. In \cite{T83}, Talenti
proved $|\nabla u|^{-1}k$ is a harmonic function. Our Corollary \ref{2-dim}
is a minimal graph version of Talenti's result. So from Corollary \ref{2-dim}
above we can also get the upper bound estimates on the curvature of the convex
level curve of 2-dimensional minimal graph with boundary data.
\end{Rem}

\section{Proof of theorem \ref{PoissonThm}}
In this section we will examine the Poisson equation in detail by setting $F^{\a\b}=\delta_{\a\b}$ and $\rho(t)=\t \log t$. Obviously, we have $\rho' u_n^{2} =\t,\ \rho'' u_n^{4}=-\t$. Recalling the relation $u_{ii}=-u_na_{ii}$, we have
\begin{align*}
u_{nn}=f+u_n\sigma_1, \qquad{\text{where}} \ \ \sigma_1=\sum_\i a_{ii}.
\end{align*}
Similar to the last section, set
$$\varphi=\t \log(|\nabla u|^2)+ \log K(x).$$
For suitable choice of $\t$ we will derive the following differential inequality
\begin{align*}
\Delta\varphi\leq 0 \qquad \mod\nabla\varphi \quad\text{in}\quad\Omega,
\end{align*}
where we modify the terms involving $\nabla \varphi$ with locally bounded coefficients. We shall complete our calculation under the normal coordinates at the fixed point $x_o$.

Recalling the following formula we have obtained in \eqref{H1234}
\begin{align}\label{H1234b}
\Delta\varphi=L_1+L_2+L_3+L_4,
\end{align}
we will treat the four terms on the right-hand side of \eqref{H1234b} respectively.

For the term $L_1$, we immediately have
\begin{align}\label{H1ok}
\begin{split}
L_{1}=&-6u_n^{-3}f\sum_\i
a^{ii}u_{ni}^2 +6u_n^{-2}\sum_\i a^{ii}u_{ni} \n_if+(2\t-n+1)u_n^{-1}\n_n f \\
&-u_n^{-1}\sum_\i a^{ii} \n_{ii}f.
\end{split}
\end{align}

For the term $L_2$, by \eqref{uiia} we have
\begin{align}\label{H2ok}
\begin{split}
L_2=&- \sum_{1\leq i,j,k\leq n-1} a^{ii}a^{jj}a_{ij, k}^2-
\sum_{1\leq i,j\leq n-1} a^{ii}a^{jj}a_{ij,
n}^2-4u_n^{-2}\sum_{1\leq i,j\leq n-1}a^{ii}u_{ni}u_{jj i}\\&
-4u_n^{-1}\sum_{\i}u_{iin}\\
=&- \sum_{1\leq i,j,k\leq n-1} a^{ii}a^{jj}a_{ij, k}^2 -4u_n^{-2}\sum_{1\leq i,j\leq n-1}a^{ii}u_{ni}(-u_na_{jj,i}+2u_n^{-1}u_{nj}u_{ji} +u_n^{-1}u_{ni}u_{jj})\\
&- \sum_{1\leq i,j\leq n-1} a^{ii}a^{jj}a_{ij, n}^2 -4u_n^{-1} \sum_{\i}(-u_na_{ii,n}+2u_n^{-1}u_{ni}^2 +u_n^{-1}u_{nn}u_{ii})\\
=&-  \sum_{1\leq i,j,k\leq n-1} a^{ii}a^{jj}a_{ij, k}^2- \sum_{1\leq
i,j\leq n-1} a^{ii}a^{jj}a_{ij, n}^2+4u_n^{-1}\sum_{1\leq i,j\leq
n-1}a^{ii}u_{ni}a_{jj, i}\\
&+4\sum_{\i}a_{ii,n}+4u_n^{-2}\sigma_1\sum_{\i}a^{ii}u_{ni}^2
+4u_n^{-1}u_{nn}\sigma_1.
\end{split}
\end{align}
Similarly, we can obtain
\begin{align}\label{H34ok}
\begin{split}
L_3=&~u_n^{-2}\bigg\{ 2\t
u_{nn}^2+(4+2\t)\sum_{\i}u_{ni}^2-6\sigma_1\sum_{\i}a^{ii}u_{ni}^2-2u_nu_{nn}\sigma_1\\
&+[2\t-(n+1)]\sum_{\i}u_{ni}^2+[2\t-(n+1)]u_n^2\sum_{\i}a_{ii}^2\bigg\},\\
L_4=&-2u_n^{-1}\sum_{\i}u_{ni}\varphi_{i}-2u_n^{-1}u_{nn}\varphi_{n}.
\end{split}
\end{align}
Combining \eqref{H1ok}--\eqref{H34ok}, we have
\begin{align*}
u_n^{2}\Delta\varphi=\overline{P}_0+\overline{P}_1+\overline{P}_2+\overline{P}_3,
\end{align*}
where
\begin{align*}
\overline{P}_0=&-6u_n^{-1}f\sum_\i a^{ii}u_{ni}^2 +6\sum_\i a^{ii}u_{ni} \n_if +(2\t-n+1)u_n\n_n f -u_n\sum_\i a^{ii} \n_{ii}f, \\
\overline{P}_1=&- u_n^{2}\sum_{1\leq i,j,k\leq n-1} a^{ii}a^{jj}a_{ij, k}^2 +4u_n\sum_{1\leq i,j \leq n-1} a^{ii}u_{ni}a_{jj, i} -u_n^{2}\sum_{1\leq i,j\leq n-1} a^{ii}a^{jj}a_{ij, n}^2 +4u_n^2\sum_{\i}a_{ii,n},\\
\overline{P}_2=&~(2\t-n-1)u_n^2\sum_{\i}a_{ii}^2-2\sigma_1\sum_{\i}a^{ii}u_{ni}^2 +(2\t+2)u_n^{2}\sigma_1 ^2 +\Big[4\t-(n-3)\Big]\sum_{\i}u_{ni}^2\\
&+2\t f^2+(4\t+2)u_nf\sigma_1,\\
\overline{P}_3=&-2u_n\sum_{\j}u_{nj}\varphi_j-2(u_nf+u_n^2\sigma_1)\varphi_n.
\end{align*}

To deal with the term $\overline{P}_1$, we  may set $A=C=D=u_n^2, B=u_n$ in
Lemma \ref{Form2}. Let us denote by $(\overline{Q}_0)$ the terms in $\overline{P}_0$; $(\overline{Q}_1)$ the terms involving $a_{ii,n}(2\leq i \leq n-1)$; $(\overline{Q}_2)$ the terms involving $a_{ii,1}(2\leq i \leq n-1)$; $(\overline{Q}_3)$ the terms involving $a_{ii,j}(2\leq i, j \leq n-1)$;  $(\overline{Q}_4)$ the terms involving $\n \varphi$ and $(\overline{Q}_5)$ all of the rest terms. More precisely, we have
\begin{align}\label{biaodashi}
\begin{split}
u_n^2\Delta\varphi=\overline{Q}_0+\overline{Q}_1+\overline{Q}_2+\overline{Q}_3 +\overline{Q}_4+\overline{Q}_5,
\end{split}
\end{align}
where
\begin{align*}
\overline{Q}_0=&-6u_n^{-1}f\sum_\i a^{ii}u_{ni}^2 +6\sum_\i a^{ii}u_{ni} \n_if+(2\t-n+1)u_n\n_n f -u_n\sum_\i a^{ii} \n_{ii}f,\\
\overline{Q}_1=&-\bigg(\sum_{2\leq i \leq n-1}u_na^{ii}a_{ii,n} \bigg)^2 +4u_n\sum_{2\leq i \leq n-1}(a_{ii}-a_{11} -\t  \sigma_1-\t u_n^{-1} f) \cdot(u_na^{ii}a_{ii,n}) \\
&-\sum_{2\leq i\leq n-1}(u_na^{ii}a_{ii,n})^2,
\end{align*}
\begin{align*}
\overline{Q}_2=&-\bigg(\sum_{2\leq i\leq n-1}u_na^{ii}a_{ii,1}\bigg)^2
+4u_{n1}\sum_{2\leq i\leq n-1}(a_{ii}a^{11}-1-\t) \cdot(u_na^{ii}a_{ii,1})\\
&-\sum_{2\leq i \leq n-1}(1+2a_{ii}a^{11})\cdot(u_na^{ii}a_{ii, 1})^2,\\
\overline{Q}_3=&\sum_{2\leq j\leq n-1}\bigg\{ -(1+2a_{11}a^{jj})\cdot
\bigg(\sum_{2\leq i \leq n-1}u_na^{ii}a_{ii,j}\bigg)^2-\sum_ {\substack{2\leq i\leq n-1\\ i\neq j}}2a_{ii}a^{jj}\cdot(u_na^{ii}a_{ii, j})^2\\
&\hspace{0.6cm}-\sum_{2\leq i\leq n-1}(u_na^{ii}a_{ii,j})^2 +4u_{nj}\sum_{2\leq i \leq n-1}[a_{ii}a^{jj}-\t-(1+2\t)a_{11}a^{jj}] \cdot(u_na^{ii}a_{ii,j}) \bigg\},\\
\end{align*}
and
\begin{align*}
\overline{Q}_4=&~2u_n^2\varphi_1\sum_\ii a^{ii}a_{ii,1}+4\t u_n u_{n1}
\varphi_1+2u_n^2\sum_{2\leq i,j\leq n-1}
(1+2a_{11}a^{jj}) a^{ii}a_{ii,j}\varphi_j\\
& +4\t u_n\sum_{2\leq j\leq n-1}(1+2a_{11}a^{jj})u_{nj}\varphi_j
+4u_na_{11}\sum_{2\leq j\leq n-1} a^{jj}u_{nj}\varphi_j -u_n^2\varphi_1^2 -u_n^2\varphi_n^2\\
&+2\Big[u_n^2\big(\sum_{\ii}a^{ii}a_{ii,n}+2\t
u_n^{-1}u_{nn}\big)+2u_n^2a_{11}\Big]\varphi_n-u_n^2\sum_{2\leq
j\leq n-1}(1+2a_{11}a^{jj}) \varphi_j^2+\overline{P}_3,\\
\overline{Q}_5=&~-u_n^2\sum_{\substack{1\leq i, j, k\leq n-1\\ i\neq j, j\neq
k, k\neq i}}a^{ii}a^{jj}a_{ij, k}^2-u_n^2\sum_{\substack{1\leq
i,j\leq n-1\\i\neq j}} a^{ii} a^{jj} a_{ij,n}^2
+(2\t-n-1)u_n^2\sum_\j a_{jj}^2\\
&+\sum_{2\leq j\leq n-1}\big(-4\t^2 +4\t-n+3 -8\t^2a_{11}a^{jj} -8\t
a_{11}a^{jj} -2\sigma_1a^{jj} \big) u_{nj}^2-8\t u_n^2a_{11}\sigma_1\\
&+\left(-4\t^2-4\t-n+3-2\sigma_1a^{11}\right) u_{n1}^2 +\big(-4\t^2
+2\t+2\big)u_n^2\sigma_1^2 +\big(-4\t^2+2\t\big)f^2\\
&+\big(-8\t^2+4\t+2\big)u_nf\sigma_1 -8\t u_n a_{11}f.
\end{align*}

In the following, we shall maximize the terms $\overline{Q}_1, \overline{Q}_2$ and $\overline{Q}_3$ via Lemma \ref{maximum} for different choice of parameters.

At first let us examine the term $\overline{Q}_1$. For $\ii$, set $X_i =u_na^{ii} a_{ii,n}$, $\lambda=1$, $\mu=u_n$, $b_i=1$ and $c_i=a_{ii} -a_{11}-\t  \sigma_1-\t u_n^{-1} f$. By Lemma \ref{maximum}, we have
\begin{align}\label{punlast}
\begin{split}
\overline{Q}_1\leq&~\frac{4(n-2)}{n-1}\t^2u_n^2\sigma_1^2 -\frac{8}{n-1}\t
u_n^2 \sigma_1^2 +8\t u_n^2a_{11}\sigma_1 -\frac{4}{n-1}u_n^2\sigma_1^2 +4u_n^2\sum_\j a_{jj}^2\\
&+\frac{4(n-2)}{n-1}\t^2f^2 +\frac{8n-16}{n-1} \t^2u_n\sigma_1f -\frac{8}{n-1}\t u_n\sigma_1f+ 8\t u_na_{11}f.
\end{split}
\end{align}

For the term $\overline{Q}_2$, set $X_i=u_na^{ii}a_{ii,1}$, $\lambda=1$, $\mu=u_{n1}$, $b_i=1+2a_{ii}a^{11}$ and $c_i=a_{ii}a^{11}-1-\t$ where $\ii$. Also by Lemma \ref{maximum} and the same discussion as before, we can get
\begin{align}\label{pun1last}
\overline{Q}_2\leq \bigg[\frac{4(n-2)}{n-1}(1+\t)^2+2\sigma_1a^{11}-2\bigg] u_{n1}^2.
\end{align}

In a similar way, for the term $\overline{Q}_3$, just copying the calculation as the minimal graph, we can derive that
\begin{align}\label{punjlast}
\overline{Q}_3\leq \sum_{2\leq j\leq n-1}\bigg(\frac{4n-8}{n-1}\t^2 -\frac{8}{n-1}\t +\frac{2n-6}{n-1} +2\sigma_1a^{jj} +8\t^2a_{11}a^{jj} +8\t a_{11}a^{jj}\bigg)u_{nj}^2.
\end{align}

Combining \eqref{biaodashi}--\eqref{punjlast}, we finally obtain
\begin{equation}\label{Last}
\begin{split}
u_n^2\Delta \varphi\leq &~\big(2\t-n+3\big)u_n^2\sum_\j a_{jj}^2
+\left[-\frac{4}{n-1}\t^2 -\frac{8}{n-1}\t +2\t -\frac{4}{n-1} +2\right]u_n^2\sigma_1^2\\
&+\sum_{1 \leq j \leq n-1}\big(-\frac{4}{n-1}\t^2 -\frac{8}{n-1}\t
+4\t-\frac{4}{n-1}-n+5 \big) u_{nj}^2\\
&+\left(-\frac{4}{n-1}\t^2+2\t\right)f^2 +\left(-\frac{8}{n-1}\t^2
-\frac{8}{n-1}\t +4\t+2 \right)u_nf\sigma_1\\
&-6u_n^{-1}f\sum_\i a^{ii}u_{ni}^2 +6\sum_\i a^{ii}u_{ni}\n_if +(2\t-n+1)u_n\n_nf\\
&-u_n\sum_\i a^{ii} \n_{ii}f\quad \mod \n \varphi,
\end{split}
\end{equation}
where we have modified the terms involving $\nabla \varphi$ with locally bounded coefficients.

To estimate the term $\overline{Q}_0$, we shall make the following choice.

{\bf Case (i):}

For $f=f(u)$, since at $x_o$,
$$\n_if =0, \quad \n_nf=u_n f_u \quad \text{and}\quad \n_{ii}f = -u_n f_u a_{ii}\quad \text{for}\quad 1\leq i\leq n-1,$$
we have
\begin{align}\label{fua}
\begin{split}
\overline{Q}_0=2\t u_n^2f_u -6u_n^{-1}f\sum_\i a^{ii}u_{ni}^2.
\end{split}
\end{align}

If $f_u\geq0$, then we choose $\t=-1$. By \eqref{Last}--\eqref{fua}, we obtain
\begin{align*}
\begin{split}
u_n^2\Delta \varphi \leq& -(n-1)\sum_\j a_{jj}^2 -(n-1)\sum_\j u_{nj}^2 -\bigg(\frac{4}{n-1}+2\bigg)f^2\\
&-2u_nf\sigma_1 -2u_n^2f_u -6u_n^{-1}f\sum_\i a^{ii}u_{ni}^2\\
\leq & 0 \quad \mod \n \varphi.
\end{split}
\end{align*}

If $f_u\leq 0$, we can set $\t=\frac {n-1}{2}$ to derive
\begin{align*}
\begin{split}
u_n^2\Delta \varphi \leq &~2u_n^2\sum_\j a_{jj}^2 -\bigg( \frac{4}{n-1} +2\bigg)u_n^2\sigma_1^2-\frac{4}{n-1}\sum_{1 \leq j \leq n-1} u_{nj}^2\\
&-2u_nf\sigma_1+(n-1) u_n^2f_u -\frac{6}{u_n}f\sum_\i a^{ii}u_{ni}^2\\
\le& 0\quad \mod \n \varphi.
\end{split}
\end{align*}

{\bf Case (ii):}

For $f=f(x)$, at the considered point $x_o$, we have
\begin{align*}
\n_if=f_{x_i},\  \n_nf=f_{x_n} \ {\text{and}} \  \n_{ii}f=f_{x_ix_i} \quad {\text{for}}1\leq i\leq n-1,
\end{align*}
hence
\begin{align*}
\overline{Q}_0=&-6u_n^{-1}f\sum_\i a^{ii}u_{ni}^2 +6\sum_\i
a^{ii}u_{ni}f_{x_i} +(2\t-n+1)u_nf_{x_n}\\
&-u_n\sum_\i a^{ii}f_{x_ix_i}.
\end{align*}

By setting $\t = \dfrac{n-1}{2}$ and applying the condition that
\begin{align*}
t^3f(x)\ {\text{is convex with respect to}} (x,t)\in\Omega\times (0,+\infty),
\end{align*}
one can easily check
\begin{align*}
6u_n^{-1}f u_{ni}^2- 6u_{ni}f_{x_i}+u_n f_{x_ix_i}\ge 0
\end{align*}
for each $1\leq i \leq n-1$.

Therefore, we have
$$\overline{Q}_0 \leq 0.$$
It follows that
\begin{align*}
\begin{split}
u_n^2\Delta \varphi \leq &~2u_n^2\sum_\j a_{jj}^2 -\bigg(\frac{4}{n-1} +2\bigg) u_n^2\sigma_1^2-\frac{4}{n-1}\sum_{1 \leq j \leq n-1} u_{nj}^2-2u_nf\sigma_1\\
\leq &0\quad \mod \n \varphi.
\end{split}
\end{align*}
The proof of the Theorem \ref{PoissonThm} is now completed. \qed

\section{Proof of the Corollary}
\setcounter{equation}{0} \setcounter{theorem}{0}

Let $\Omega_0$ and $\Omega_1$ be bounded smooth convex domains in $\R^n, n \ge 2$, $o\in \bar\Omega_1\subset\Omega_0$. Let $u$ satisfy
\begin{equation}\label{quasilinear}
\left\{
\begin{array}{lcl}
              \Delta u = f(u)   &\text{in}&  \Omega=\Omega_0\backslash\bar\Omega_1, \\
                     u = 0   &\text{on}&  \partial \Omega_0,\\
                     u = 1   &\text{on}&  \partial \Omega_1,
\end{array} \right.
\end{equation}
where $f\in C^2([0, 1])$ is a nonnegative and non-decreasing function with $f(0)=0$.

 From \cite{CS82} we know that $|\nabla u|\neq 0$ in
$\Omega$ if $u$ is a solution of \eqref{quasilinear}. In the following, we
shall prove a lemma on the monotonicity of the norm of the gradient along the
gradient direction, which also appeared in \cite{LS}. Using this observation,
we prove the Corollary~\ref{CorCaffarelli}.

\begin{Lem}\label{IncreaseGradident}
Let $u$ satisfy \eqref{quasilinear}. Then $|\nabla u|$ strictly increases in the direction $\nabla u$. It follows that $|\nabla u|$ attains its minimum on $\partial \Omega_0$, and attains its maximum on $\partial \Omega_1$.
\end{Lem}
\begin{proof}
By the Caffarelli-Spruck's \cite{CS82}, the level sets of $u$ are strictly convex with respect to the normal direction $\nabla u$. At any fixed point $x_o\in\Omega$, we may let $u_i=0 (1\leq i\leq n-1)$ and $u_n=|\nabla u|>0$ by rotation. Let $H$ be the mean curvature of the level sets with respect to the normal direction $\nabla u$. Then $\eqref{quasilinear}$ implies
\begin{align*}
u_{nn}=-\sum_\i u_{ii}+f=u_nH+f,
\end{align*}
hence
\begin{align*}
\sum_{1\leq \a\leq n}(|\nabla u|^2)_\a u_\a =2u_n^2u_{nn} =2u_n^2(u_nH+f)>0,
\end{align*}
where the last inequality is due to the strict convexity of the level sets.
\end{proof}

Now we give the proof of Corollary~\ref{CorCaffarelli}.
\begin{proof}
If $u$ is the smooth solution of $\eqref{quasilinear}$, then from the
Caffarelli-Spruck's \cite{CS82} we know the level sets of $u$ are strictly convex with respect to normal direction $\n u$. From the Theorem \ref{PoissonThm}, the function $$|\n u|^{-2}K$$ attains it minimum on the boundary. Since $|\nabla u|$ attains its minimum on $\partial \Omega_0$, and attains its maximum on $\partial \Omega_1$, we have the following estimate
\begin{align*}
\min_{\Omega} K \geq \bigg( \frac{\min_{\partial \Omega_0}|\n u|} {\max_{\partial\Omega_1}|\n u|}\bigg)^2 \min_{\partial \Omega}K.
\end{align*}
Then we complete the proof of the Corollary~\ref{CorCaffarelli}.
\end{proof}

\end{document}